\documentclass[a4paper,reqno,12pt]{amsart}
\usepackage{rr_prefs}
\author[Raj Rao Nadakuditi]{Raj Rao Nadakuditi}\address{Raj Rao Nadakuditi, Department of Electrical Engineering and Computer Science, University of Michigan, 1301 Beal Avenue, Ann Arbor, MI 48109. USA.}
\email{rajnrao@eecs.umich.edu}
\urladdr{http://www.eecs.umich.edu/\~\/rajnrao/}

\thanks{This work was supported by an ONR Young Investigator Award  N000141110660, an AFRL subcontract from Solid State Scientific (PM Mike Noyola), AFOSR Young Investigator Award FA9550-12-1-0266 and a ARO MURI grant W911NF-11-1-0391. The author thanks Jeff Fessler for teaching him about the role of non-convex penalty functions and Iain Johnstone for suggesting the phrase `observable solution to an unobservable problem' to emphasize the surprising nature of the optimum. We thank Florent Benaych-Georges for many stimulating conversations and Jack Silverstein for discussions on the importance of the low-coherence condition in the missing entries setting. W e thank Brendan Farrell for his suggestion to analyze the missing data problem using the signal-plus-noise-plus-small-perturbation framework and the reviewers for their suggestions. A short (conference) version containing some of the ideas in this paper first appeared in \cite{nadakuditi2011exploiting}.}

\title[Low rank matrix approximation]{OptShrink: An algorithm for improved low-rank signal matrix denoising by optimal, data-driven singular value shrinkage}
 %rand mat

\begin{document}

\begin{abstract}
The truncated singular value decomposition (SVD) of the measurement matrix is the optimal solution to the \textit{representation} problem of how to best approximate a noisy measurement matrix using a low-rank matrix. Here, we consider the (unobservable) \textit{denoising} problem of how to best approximate a low-rank signal matrix buried in noise by optimal (re)weighting of the singular vectors of the measurement matrix. We exploit recent results from random matrix theory to exactly characterize the large matrix limit of the optimal weighting coefficients and show that they can be computed directly from data for a large class of noise models that includes the i.i.d. Gaussian noise case.

Our analysis brings into sharp focus the shrinkage-and-thresholding form of the optimal weights, the non-convex nature of the associated shrinkage function (on the singular values) and explains why matrix regularization via singular value thresholding with convex penalty functions (such as the nuclear norm) will always be suboptimal. We validate our theoretical predictions with numerical simulations, develop an implementable algorithm (\textrm{OptShrink}) that realizes the predicted performance gains and show how our methods can be used to improve estimation in the setting where the measured matrix has missing entries.
\end{abstract}

%\tableofcontents
\maketitle

% !TeX root = lowrankapprox.tex
\section{Introduction}
Techniques for low-rank signal matrix extraction from a signal-plus-noise matrix appear prominently in many statistical signal processing \cite{tufts1993estimation,scharf1991svd,jolliffe2002principal}, machine learning \cite{drineas2006fast,kannan2009spectral}, estimation and classification applications \cite{klema1980singular}. In many applications, the low-rank approximation is the first step in an inferential process  (see, for e.g. \cite{tipping1999mixtures,dasgupta1999learning,sanjeev2001learning,vempala2002spectral,kannan2005spectral,hsu2012learning}). These techniques are necessary whenever the $n \times m$ signal-plus-noise data or measurement matrix formed by, for example lining up the $m$ samples or measurements of $n \times 1$ observation vectors alongside each other, can be modeled as
\begin{equation}\label{eq:model}
\widetilde{X} = \sum_{i=1}^{r} \theta_{i} u_{i} v_{i}^{H} + X,
\end{equation}
where $^H$ denotes the conjugate transpose and $u_{i}$ and $v_{i}$ are left and right ``signal'' singular vectors associated with singular values $\theta_i$ of the signal matrix \be \label{eq:signal model} S =  \sum_{i=1}^{r} \theta_{i} u_{i} v_{i}^{H}\ee  and $X$ is the noise-only matrix of random (not necessarily i.i.d.) noises. These models also arise in other graph signal processing type settings; see for example \cite[Text before (9)]{nadakuditi2012graph},  \cite[Section V]{nadakuditi2013spectra}, \cite[Section III.A]{karger2011budget} or the various models described in \cite{chatterjee2012matrix}.

Relative to this model the objective is to form an estimate of the low-rank signal matrix assuming, for now, that its rank $r$ is known. The truncated singular value decomposition (SVD) plays a prominent role in a widely-used `optimal' solution to a problem that is addressed by the famous Eckart-Young-Mirsky (henceforth, EYM) theorems \cite{eckart1936approximation,mirsky1960symmetric,golub1987generalization}. Specifically, if  $\fronorm{\cdot}$ denotes the matrix Frobenius norm then the solution to the constrained optimization problem
\be \label{eq:eym formulation}
 \widehat{S}_{\rm eym} = \argmin_{\rank(S) = r} \fronorm{\widetilde{X} - S},
\ee
is given by
$$\widehat{S}_{\rm eym} = \sum_{i=1}^{r} \widehat{\sigma}_{i}\widehat{u}_i \widehat{v}_i^{H},$$ where $\widetilde{X} = \sum_i \widehat{\sigma}_i \uhat{i} \vhat{i}^H$ is the SVD of $\widetilde{X}$. This is also the maximum likelihood (ML), rank $r$ estimate when $X$ is assumed to be a matrix with i.i.d. Gaussian entries since the negative log-likelihood function is precisely the right hand side of (\ref{eq:eym formulation}). Its use is also justified in the small $n$, large $m$ (or vice versa) regime, whenever local asymptotic normality \cite{le1960locally} has `kicked in'.

A natural extension is to consider settings where the signal matrix is low rank and has some additional exploitable structure. Examples include low-rank and sparse (see the body of work on sparse principal component analysis. e.g.  \cite{zhang2002low,zou2006sparse,johnstone2009consistency,d2008optimal,jenatton2010structured,rohde2011estimation,zhang2012sparse,birnbaum2012minimax}),  low rank and Toeplitz structured (see e.g. \cite{tufts1993estimation,cadzow1991enhanced,wikes1988iterated}), low rank and Hankel structured \cite{li1997parameter} and low rank and nonnegative \cite{boutsidis2008svd,langville2006initializations,chu2004optimality}; see \cite{chu2003structured,markovsky2008structured} for an excellent overview of these methods and additional references. As expected,  by exploiting structure in the signal matrix we can improve estimation performance relative to the EYM estimator which assumes no structure besides the low-rank condition.

\subsection{Denoising by optimally weighted approximation}
Here we place ourselves in the setting where \textit{no structure is assumed in the low-rank signal matrix} and ask how the EYM estimator can be improved. The starting point for our investigation is the observation that as formulated in (\ref{eq:eym formulation}), the EYM estimator solves the \textit{representation problem} of finding the best rank $r$ approximation of the \textit{signal-plus-noise} measurement matrix. It says nothing about the \textit{denoising} problem of how to best estimate the low-rank \textit{signal matrix}, even though practitioners sometimes invoke it as though it does. Thus we should not expect the EYM estimator to be the optimal solution to the denoising problem.

Let $||w||_{\ell_0}=|\{\#i: w_i \neq 0\}|$ so that  $||w||_{\ell_0} =r$  denotes a vector $w$ with $r$ non-zero entries. In this paper, we consider variations of the denoising problem  formulated as a weighted approximation problem of the form
\begin{equation}\label{eq:wopt generic}
\wopt{} :=  \argmin_{||w||_{\ell_0} = r } \fronorm{\sum_{i=1}^{r} \theta_i u_i v_{i}^{H} - \sum_{i} w_{i} \widehat{u}_{i} \widehat{v}_{i}^{H}}.\ee
Note that in (\ref{eq:wopt generic}), we are trying to approximate the unknown signal matrix using the singular vectors estimated from the noisy measurement matrix. Our setup is different from other weighted low-rank approximation problems considered in the literature as in \cite{srebro2003weighted}, which involve weighted modifications of the problem in (\ref{eq:eym formulation}). In our formulation, setting $w_i = \widehat{\sigma}_i$ recovers the EYM estimator so that by inspecting the solution we can directly assess when and the extent to which the EYM estimator will be suboptimal.

We prove, using recent results from random matrix theory \cite{benaych2011svd}, that for a large class of noise models, which includes but goes well beyond the i.i.d. Gaussian model, we can compute $\wopt{}$ in closed-form in the large matrix limit. The computation shows that $\wopt{}$ depends only on (an integral transform of) the limiting singular value distribution of the noise-only matrix $X$. We then exploit this fact to develop a concrete algorithm for computing a consistent (in a sense we make precise)  estimate of the limiting oracle solution directly from measurement matrix.
\subsection{Form of the optimal shrinkage-and-thresholding operator}
The analysis shows that $\wopt{i}$ takes the form of a shrinkage-and-thresholding operator (on the singular values of $\wtX$) that is completely characterized by the limiting singular value distribution of the noise-only  matrix. The resulting shrinkage function is non-convex with $\wopt{i} \approx \widehat{\sigma}_{i}(1-O(1/\widehat{\sigma}_{i}^2))$ for large $\widehat{\sigma}_{i}$ and $\wopt{i} \to 0$ for $\widehat{\sigma}_{i} \leq b + o(1)$  where $b$ is a critical threshold that depends on the limiting noise-only singular value distribution.

The shrinkage portion of the solution arises because $\widehat{\sigma}_{i}$ is positively biased relative to $\theta_i$ and because the corresponding singular vectors of $\wtX$ are biased, noisy estimates of the (true) singular vectors of the latent signal matrix \cite{benaych2011svd}. The thresholding portion of the solution arises because of a phase transition in the `informativeness' of the estimated singular vectors, relative to the latent singular vectors whereby for $\theta_i > \theta_c$ inner-products of the form $\uhatu{i}{i}$ and $(\vhatv{i}{i})$ are $O(1)$ and tend to a constant, while for  $\theta_i < \theta_c$, inner-products of the form $\uhatu{i}{i}$ and $(\vhatv{i}{i})$ are $o(1)$  and tend to zero.

Our analysis of the structure of the optimal solution 1) brings into sharp focus  the form of the optimal  shrinkage-and-thresholding operator, 2) provides insight on why the EYM estimator is near optimal in the low noise regime but sub-optimal in the moderate to high noise regime and 3) explains why we can expect that soft thresholding (of singular value) operators with convex penalty functions (such as the nuclear norm \cite{cai2010singular}) that are tuned to be near-optimal in the small $\theta_i$ regime will be suboptimal in the large $\theta_i$ regime (and vice versa).

%  shrinkage  regularization \cite{cai2010singular,negahban2011estimation,kakade2012regularization,negahban2012restricted}

\subsection{Mitigating the effect of rank over-estimation}\label{sec:rank est}

It is a delightful fact that even though the optimization problem in (\ref{eq:wopt generic}) is unobservable, because it depends on the unknown matrix we are trying to estimate, the optimal solution itself is computable. We assume no structure, other than low rank, on the signal matrix; the exploitable structure is present in the `noise portion' of the eigen-spectrum, \ie, the $\min(m,n)-r$ singular values of $\wtX$.

This makes contact with the important question of how to estimate $r$ in (\ref{eq:model}) so that one may distinguish the `signal portion' of the eigen-spectrum from the `noise portion'. The problem has been completely solved for the setting where $X$ has i.i.d. Gaussian entries. In this setting, the recentering and rescaling constants that must be applied to the largest eigenvalue of $XX^H$  to produce the Tracy-Widom distribution can be precisely characterized and used to set the appropriate threshold; see \cite{combettes92,johnstone2001distribution,baik2005phase,baik2006eigenvalues,johnstone2006high,el2007tracy,paul2007asymptotics,ulfarsson08,nadakuditi2008sample,onatski09,onatski2010determining,kritchman2008determining,kritchman2009non,nadler2010nonparametric,passemier12}. Recent work on the universality of this limiting distribution
\cite{soshnikov2002note,feral2009largest,erdHos2012universality,pillai2011universality,cacciapuoti2012local,peche2009universality} provides a rigorous justification for using  essentially the same method  in the non-Gaussian setting.

Similarly, when the columns of $X$ are i.i.d. and each column has a (non-identity) population covariance matrix with a known (limiting) eigen-distribution, then the results in \cite{el2007tracy} facilitate computation of the appropriate threshold for distinguishing the `noise portion' of the eigen-spectrum from the `signal portion'.

If the form of population covariance matrix is misspecified then applying the tests based on this theory will lead to an overestimation of the rank of the signal matrix. Developing robust estimators of the signal rank that ``work'' without having to specify the symmetry structure (e.g. i.i.d. elements, i.i.d. columns, variance profile, etc.) of the noise random matrix remains an important open problem. Such estimators will have to exploit (symmetry-independent) `universal' features of the spectrum in a way that present estimators do not.

This is where the algorithm we have developed really shines. Our algorithm takes as its input an estimate of the rank of the signal matrix and returns a (re)weighted approximation that largely mitigates the effect of rank overestimation in a manner that the EYM estimate cannot. Thus, advances in robust rank estimation when used with our algorithm will lead to improved signal matrix approximation. If the rank is correctly estimated, then the algorithm will better estimate weak subspace components of the signal matrix than the EYM algorithm.

\subsection{Contributions}
Characterizing the limiting solution of (\ref{eq:wopt generic}), computing the resulting limiting squared error, quantifying the improvement relative to the EYM estimator and developing an implementable algorithm that realizes these performance gains are the main contributions of this paper. Some of the ideas in this paper were initially presented in a conference paper by the author \cite{nadakuditi2011exploiting}, in the context of the i.i.d Gaussian noise setting. This version goes beyond the Gaussian setting considered there. We also treat the setting where measurement matrix has missing entries, as considered in \cite{chen2004recovering,fazel2008compressed,candes2010matrix,candes2010power,candes2009exact,keshavan2010matrix}. In addition to rigorous results, we formulate some (empirically validated and theoretically justified) conjectures for the structure of the solution for various `rank-regularized' variations of (\ref{eq:wopt generic}).

In related work, Hachem et al \cite{hachem2012subspace} looked at the problem of structured subspace estimation arising in the context of parameter estimation in large arrays. They propose an oracle solution \cite[Equation (13), pp. 435]{hachem2012subspace} and analyze its first and second order performance in the context of the MUSIC direction-of-arrival estimator.

If we were to apply the ideas and techniques developed in this paper to the problem
$$\wopt{} :=  \argmin_{||w||_{\ell_0} = r }|| \sum_{i=1}^{r} u_{i} u_i^{H} - w_i \uhat{i} \uhat{i}^H ||_{F}^{2},$$
then, we would recover a solution that corresponds to their oracle solution. Here, we consider the problem of estimating the low-rank matrix; our results and our new algorithm can be analyzed using the techniques in \cite{hachem2012subspace} to provide insights on the first and second order convergence properties. We leave the extension of our techniques to the estimation of projection matrices is relatively straightforward as an exercise to the reader.

The paper is organized as follows. The setup, the main theoretical results and a new algorithm based on the theoretical analysis are presented in Section \ref{sec:main results}.  Simulation results to validate the theoretical predictions and a comparison of our method to other matrix regularization methods are contained in Section \ref{sec:discussion}.

\section{Main results and a new algorithm}\label{sec:main results}

\subsection{Setup and Notation}

Let  $X_n$ be  an $n \times m$ ($n \leq m$, without loss of generality\footnote{We choose this convention to simplify the definition of the empirical singular value distribution.}) random matrix whose ordered singular values we denote by $\si_{1}(X_n) \geq \cdots \geq \si_{n}(X_n)$. Let $\mu_{X_{n}}$ be the empirical singular value distribution, \ie, the probability measure defined as
\[
\mu_{X_{n}} = \frac{1}{n} \sum_{i=1}^{n} \delta_{\si_{i}(X_{n})}.
\]
Assume that the probability measure $\mu_{X_{n}}$ converges almost surely weakly, as $n,m \longrightarrow \infty$, to a non-random compactly supported  \pro measure $\mu_X$ that is supported on $[a,b]$. We assume that  $\si_{1} \convas b$, where $\convas$ denotes almost sure convergence. These conditions are satisfied by the model where $X_n$ has i.i.d. entries mean zero entries with variance $1/m$ and bounded higher order moments.

For a given $r \ge 1$, let $\tta_1 > \cdots >  \tta_r >0$ be deterministic non-zero real numbers, chosen independently of $n$.  For every $n$, let $S_n$ be an $n \times m$ signal matrix   having rank $r$ with its $r$ non-zero distinct singular values equal to $\theta_{1}, \ldots, \theta_{r}$.

We suppose that $X_n$ and $S_n$ are independent and that $X_n$, the noise-only matrix  is bi-unitarily invariant while the low-rank signal matrix $S_n$ is deterministic.  Recall that a random matrix is said to be {\it bi-orthogonally invariant} (or {\it bi-unitarily invariant}) if its distribution is invariant under multiplication on the left and right by orthogonal (or unitary) matrices. Alternately, if $S_n$ has isotropically random right (or left) singular vectors, then $X_n$ need not be unitarily invariant under multiplication on the right (or left, resp.) by orthogonal or unitary matrices. Equivalently, $X_n$ can have deterministic right and left singular vectors while $S_n$ can have isotropically random left and right singular vectors and we would get the same result stated shortly.

A matrix $X_n$ with i.i.d. Gaussian entries satisfies these assumption; our results extend well beyond the Gaussian setting. The main advantage of modeling the noise matrices as having isotropically random singular vectors is that it allows us to characterize the solution in terms of just the (marginal) singular  value distribution of the noise-only matrix instead of having to model  the full joint distribution of the elements of the noise-only matrix.

Since the singular value distribution of the noise-only part can be estimated from the singular value distribution of the signal-plus-noise matrix, we can develop a concrete, data-driven algorithm, presented in Section \ref{sec:algorithm}, that can applied to real-world datasets to improve low-rank signal matrix recovery.

We observe a signal-plus-noise matrix $\widetilde{X}_{n}$ modeled as, 
$$\widetilde{X}_{n} = S_{n} + X_{n},$$
where the signal matrix $S$ is modeled as in (\ref{eq:signal model}).  For $i = 1, \ldots, q = \min(n,m)$, let $\uhat{i}$ and $\vhat{i}$ denote the left and right singular vectors of $\wtX$ (we suppress the subscript $n$) associated with the singular value $\widehat{\sigma}_{i}$. The solution to the  optimization problem
\be \label{eq:weym rep}\weym{} =\argmin_{||w||_{\ell_0}=r} \fronorm{\wtX - \sum_{i} w_{i} \widehat{u}_{i} \widehat{v}_{i}^{H}},\ee
is given by $\weym{i} = \widehat{\sigma}_{i}$ for $i = 1, \ldots,r$. This yields the rank $r$ signal matrix estimate $\sum_{i=1}^{r} \weym{i} \widehat{u}_{i} \vhat{i}^H$ which, by the EYM theorem, is also the solution to the representation problem in (\ref{eq:eym formulation}).

 For $w \in \mathbb{R}^l$, define the squared error as
\be \label{eq:mse w}
\mse{w} =  \fronormsq{S  - \sum_{i=1}^{l} w_{i} \,\uhat{i}\vhat{i}^{H}}.
\ee
Consider the denoising optimization problem
\be \label{eq:wopt}
\wopt{} := \argmin_{w = [w_1 ~~ \cdots w_{r}]^T \in \mathbb{R}^{r}_{+}}  \mse{w}.
\ee
We now characterize $\wopt{}$ exactly (for every $n$) and provide an expression for its limiting value. In what follows, for a function $f$ and $c\in \R$, we set $$f(c^+):=\lim_{z\downarrow c}f(z).$$

\subsection{Theoretical results}

\begin{Th}[Weighting coefficients]\label{th:wyem vs wopt}
The solution to (\ref{eq:wopt}) exhibits the following  behavior in the asymptotic regime where $n,m\to \infty$ and $n/m \to c \in [0,\infty)$. We have that for every $1\le  i \leq r$,
\flushleft a)
\be \label{eq:wopt thm}
\wopt{i} =  \left(\Re\{ \sum_{j=1}^{r} \theta_j \uhatu{i}{j} \vvhat{j}{i} \}\right)_{+} \convas   - 2 \dfrac{D_{\mu_X}(\rho_i)}{D'_{\mu_X}(\rho_i)}  \qquad {\rm if } ~ \theta_{i}^{2} > 1/D_{\mu_X}(b^+),
\ee
where $x_{+} = \max(0,x)$ and $\rho_i =  D_{\mu_X}^{-1}(1/\theta_i^2)$.
\flushleft b)
\be \label{eq:weym thm}
\weym{i} = \widehat{\sigma}_{i} \convas \begin{cases} D_{\mu_X}^{-1}(1/\theta_i^2) = \rho_i & {\rm  if} ~\theta_i^2  > 1/ D_{\mu_X}(b^+) ,\\ \\
b & {\rm otherwise.}
\end{cases}
\ee

In a) and b),  $D_{\mu_X}(\cdot)$ is the $D$-transform of $\mu_X$ defined as
$$ D_{\mu_X}(z) :=\left[\int \frac{z}{z^2-t^{2}} \ud \mu_X(t)\right] \times \left[c\int \frac{z}{z^2-t^{2}} \ud \mu_X(t)+\frac{1-c}{z}\right] \qquad {\rm for }~ z \notin  \supp \mu_X,$$
and   $D_{\mu_X}^{-1}(\cdot)$ denotes its functional inverse.

\flushleft c) A straightforward consequence of a) and b) is that 
$$\weym{i} > \wopt{i},$$
almost surely. Also,  $\weym{i} \convas \wopt{i}$ as $\theta_i \to \infty$.
\end{Th}

The emergence of the $D$ transform in the limit characterization of the EYM and optimal coefficients follows from the results in \cite{benaych2011svd}. There it was shown that, in the large matrix limit, the principal singular values and singular vectors of $\widetilde{X}$ can be completely characterized in terms of the singular values of the signal matrix and the $D$-transform of the limiting noise-only singular value distribution. This is why, in Theorem \ref{th:wyem vs wopt}, the limiting values of $\weym{i}$ and $\wopt{i}$ only depend on the singular values $\theta_{i}$ (or $\rho_{i}$) of the signal matrix and the limiting noise-only singular value distribution $\mu_X$.

The $D$-transform is the analog of the log-Fourier transform in the sense that it describes how the distribution of the singular values of the sums of `freely' independent matrices are related to the distribution of the singular values of the individual matrices \cite{b09}. In that sense it is an \textit{asymptotically sufficient statistic} and hence its appearance in Theorem \ref{th:wyem vs wopt} is rather natural. See Section 2.5 of \cite{benaych2011svd} for additional remarks.

We now characterize the limiting squared error for the optimal, EYM and other estimators with arbitrary weights.
\begin{Th}[Limiting squared error]\label{th:limiting mse}
Assuming that for $i = 1, \ldots, r$, $\theta_{i}^{2} > 1/D_{\mu_X}(b^+)$. Then in the asymptotic regime considered in Theorem \ref{th:wyem vs wopt}, the squared error, defined as in (\ref{eq:mse w}), exhibits the following limiting behavior:
$$\mse{w} \convas  \sum_{i=1}^{r} \left( \theta_i^2 + w_i^2  +  \dfrac{4 w_i }{ \theta_{i}^{2} D'_{\mu_X}(\rho_i)} \right) $$

\noindent Consequently,
\flushleft a) $$\mse{\wopt{}} \convas  \sum_{i=1}^{r} \left( \theta_i^{2}-  \dfrac{4}{(\theta_i^2 D'_{\mu_X}(\rho_i) )^{2}} \right),$$
\flushleft b) $$\mse{\weym{}} \convas  \sum_{i=1}^{r} \left( \theta_i^{2} + \rho_i^2  +  \dfrac{4 \rho_i }{ \theta_{i}^{2} D'_{\mu_X}(\rho_i)} \right).$$
More generally
\flushleft c) $$ \mse{w} - \mse{\wopt{}}  \convas \sum_{i=1}^{r} \left(w_i + \dfrac{2}{\theta_i^2 D_{\mu_X}'(\rho_i)}\right)^2$$
so that by construction
 $$ \mse{\wopt{}} < \mse{\weym{}},$$
almost surely.
\end{Th}

Theorem \ref{th:limiting mse} reveals that whenever $w_i -\wopt{i}$ is large, we can expect a significant increase in SE relative to the optimal estimator.  The next result reveals the shrinkage-and-thresholding form of the optimal estimator.\\

\begin{Th}[Shrinkage-and-thresholding form and resulting SE]\label{th:r1 soft thresholding}When $r=1$, let the  sole non-zero singular value of $S$  be denoted by $\theta$ and assume that $ D_{\mu_X}'(b^+)=-\infty$. Then  in the asymptotic regime considered, we have that
$$\wopt{1} \convas \begin{cases}
 \dfrac{-2 }{\theta^{2} D'_{\mu_X}(\rho) } & { \rm if } ~ \theta^{2} > 1/D_{\mu_X}(b^+)\\ \\
0 & {\rm otherwise,}
\end{cases}$$
where  $\rho = D_{\mu_X}^{-1}(1/\theta^2)$.

\flushleft Consequently,
$$\mse{\wopt{}} \convas  \begin{cases}
 \theta^{2}-  \dfrac{4}{(\theta^2 D'_{\mu_X}(\rho) )^{2}}  & {\rm if }~ \theta^2 > 1/D_{\mu_X}(b^+) \\ \\
\theta^2 & {\rm otherwise.}
\end{cases}
$$
whereas
$$\mse{\weym{}} \convas  \begin{cases}
\theta^{2} + \rho^2  +  \dfrac{4 \rho }{ \theta^{2} D'_{\mu_X}(\rho)}  & {\rm if }~ \theta^2 > 1/D_{\mu_X}(b^+) \\ \\
\theta^2 + b^2 & {\rm otherwise.}
\end{cases}
$$
\end{Th}

\fl Theorem \ref{th:r1 soft thresholding} shows that when $b$ (the a.s. limit of the largest noise-only singular value) is $O(1)$, we can expect an $O(1)$ decrease in SE, relative to the EYM estimator, by thresholding whenever $\theta^2 < 1/D_{\mu_X}(b^+)$. Note that when $X$ is i.i.d. Gaussian with mean zero, variance $1/m$ entries, then $b = (1+ \sqrt{c})$ and $D'_{\mu_X}(b^+)= -\infty$ so that these results apply. More generally, whenever $\mu_X$ exhibits a square-root decay at $b$ then $D'_{\mu_X}(b^+)= -\infty$  will be satisfied. Silverstein and Choi \cite{silverstein1995analysis} show that a large class of (non i.i.d.) Gaussian noise models will satisfy this condition.

\subsection{The missing data with i.i.d. noise setting}\label{sec:missing data}

We now consider the setting where $\wtX$ has missing entries so that the signal-plus-noise matrix is modeled as
\begin{equation}\label{eq:missing model}
\widetilde{X} = \left(\sum_{i=1}^{r} \theta_{i} u_{i} v_{i}^{H} + X \right) \odot M
\end{equation}
where
$$M_{ij} = \begin{cases} 1 & \textrm{ with probability } p \\
0 & \textrm{ with probability } 1- p
\end{cases}
$$
and $\odot$ denotes the Hadamard or element-wise product. Consider the optimization problem
\be \label{eq:wopt r missing}
\wopt{} := \argmin_{w = [w_1 ~~ \cdots w_{r}]^T \in \mathbb{R}^{r}_{+}}  \fronormsq{\sum_{i=1}^{r} p \, \theta_i u_i v_i^H  - \sum_{i=1}^{r} w_{i} \uhat{i}\vhat{i}^{H}}.
\ee
Note that here we are approximating $pS$ instead of $S$ as in (\ref{eq:mse w}) (so that we can use the data-driven algorithm as-is). Setting $\wopt{i} \mapsto \wopt{i}/p$ will yield a solution to the denoising problem in (\ref{eq:mse w}). Let $||w||_{\infty}= \max_{i} | w_{i} |$ denote the element of the vector $w$ with the maximum absolute value.

\begin{Th} \label{conj:missing} Assume that the singular vectors $u_{i}$ and $v_{i}$ in (\ref{eq:missing model}) satisfy a `low-coherence' condition in the following sense: we suppose that there exist non-negative constants $\eta_{u}$, $C_{u}$, $\eta_{v}$ and $C_{v}$, independent of $n$, such that  for $i = 1, \ldots, r$
\begin{equation}
 \max_{i} || u_{i} ||_{\infty} \leq \eta_{u} \dfrac{\log^{C_u} n}{\sqrt {n}} \qquad \textrm{ and } \qquad \max_{i} || v_{i} ||_{\infty} \leq \eta_{v} \dfrac{\log^{C_v} m}{\sqrt {m}}.
 \end{equation}
Let the elements of $X_{ij}$ be i.i.d with mean zero, variance $1/m$ and bounded higher order moments. Then the solution to (\ref{eq:wopt r missing}) exhibits the following limiting behavior. We have that for $p \in (0,1]$ and $i = 1, \ldots, r$  \\
\flushleft a)
$$\weym{i} = \sigma_{i}(\wtX) \convas
\begin{cases}
\sqrt{p} \cdot \sqrt{\dfrac{(1+p\,\theta_i^2)(c+p\,\theta_i^2)}{p\,\theta_{i}^{2}}}  & {\rm if }~ \theta_{i} > \dfrac{c^{1/4}}{\sqrt{p}}, \\ \\
\sqrt{p} \,(1+\sqrt{c}) & {\rm  otherwise}.\\
\end{cases}
$$
b)

$$\wopt{i}  \convas  p\, \theta_i \cdot  \sqrt{1 -\dfrac{c(1+p\,\theta_i^2)}{p\, \theta_{i}^{2} (p\,\theta_i^2+c)}} \sqrt{1 - \dfrac{c +p\,\theta_i^{2}}{p\,\theta_i^{2}(p\,\theta_i^{2} +1 )}}   \qquad  {\rm if } ~\theta_i > \dfrac{c^{1/4}}{\sqrt{p}}.$$\\

c)  When $r = 1$

$$ \wopt{i} \convas 0  \qquad {\rm if }~  \theta_{i} \leq \dfrac{c^{1/4}}{\sqrt{p}},$$
\end{Th}

%Limit it goes to $1+c$.

Theorem \ref{conj:missing} is a statement about the optimality of the shrinkage-and-thresholding form when there are missing entries in the signal-plus-noise matrix. Note that in this case, the equivalent noise-only matrix will not bi-unitarily invariant when $X$ is non-Gaussian. The proof (see Section \ref{proof:missing}), however, reveals that it asymptotically behaves as though it does so that the results of Theorems 2.1 and 2.3 still apply. Note that as a consequence, Theorem \ref{th:limiting mse} can applied to compute the result asymptotic squared error. After the submission of this paper, we learned of recent work by Shabalin and Nobel for the $p=1$ setting of Theorem \ref{conj:missing} with i.i.d. Gaussian noise; see \cite{shabalin2013reconstruction}.

%Here too, we expect significant gains in performance relative to the EYM estimator whenever $\reff < r$.

\subsection{The asymptotic equivalence of various rank-regularized estimators}

Let us define the \textit{effective rank}, $\reff$,  of the signal matrix as
\be \label{eq:reff}
 \reff = \textrm{the number of } i \in \{1, \ldots r\} \textrm{ such that } \theta_{i}^{2} > 1/ D_{\mu_X}(b^+).
\ee
Thus, the effective rank quantifies the number of singular values in the signal-plus-noise matrix $\wtX$ that are `informative', \ie, reveal the existence of a low-rank signal matrix. Clearly, $\reff \leq r$ but $\reff < r$ whenever the number of singular values that separate from the right edge $b$ of the spectrum is less than the latent signal matrix rank $r$.  The following conjecture formalizes their relation to the number of `informative' singular vectors in the signal-plus-noise matrix.

\begin{conj}[Uninformativeness below phase transition] \label{conj:partial delocalization}
Assume that $ D_{\mu_X}'(b^+)=-\infty$ and \footnote{Note that these conditions are met when $X$ has i.i.d. entries of variance $1/m$. See Theorem 2.10 of \cite{BEKYY13}.}  that for fixed  $\rhat$,
$$\max_{i} \left( \sigma_i(X) - \sigma_{i+1}(X) \right) \leq O\left( \dfrac{\log n ~{\rm factors} }{n^{2/3}}\right),$$
with very high probability. Then we have that for $\reff < i \leq \rhat$  and $j = 1, \ldots r$,
$$\max_{i,j} |\uhatu{i}{j}| \leq  O\left( \dfrac{\log n ~{\rm factors}}{n^{1/6}}\right)   \qquad ~{\rm and } \qquad  \max _{i,j}| \vvhat{j}{i}|  \leq O\left( \dfrac{\log m ~{\rm factors}}{m^{1/6}}\right) ,$$
with high enough probability that we can establish their almost sure convergence to zero.
\end{conj}
We now consider the principal rank-regularized optimization problem
\be \label{eq:wopt r}
\wopt{}(\widehat{r}) := \argmin_{w = [w_1 ~~ \cdots w_{\widehat{r}}]^T \in \mathbb{R}^{\rhat}_{+}}  \fronormsq{\sum_{i=1}^{r} \theta_i u_i v_i^H  - \sum_{i=1}^{\widehat{r}} w_{i} \uhat{i}\vhat{i}^{H}}.
\ee
We characterize the structure of the optimal estimator and the resulting MSE next.
\begin{cor}[Performance with estimated principal component rank]\label{th:perf reff}
Let $\widehat{r}$ be a fixed (with $n$) estimate of $\reff$ and $\reff$ be defined as in (\ref{eq:reff}). Then,  in the  asymptotic regime considered, assuming Conjecture \ref{conj:partial delocalization} holds, we have that

$$\wopt{i}(\rhat)  \convas
\begin{cases}
-2 \dfrac{D_{\mu_X}(\rho_i)}{D'_{\mu_X}(\rho_i)}  & {\rm for }~ i \leq \reff \\ \\
0 & {\rm otherwise.}
\end{cases}$$
and hence
$$\mse{\wopt{}} \convas  \sum_{i=1}^{\min(\reff,\rhat)} \left( \theta_i^{2}-  \dfrac{4}{(\theta_i^2 D'_{\mu_X}(\rho_i) )^{2}} \right) + \sum_{i = \min(\reff,\rhat) +1}^{\max(r,\rhat)} \theta_{i}^2,$$
whereas
$$\mse{\weym{}} \convas  \sum_{i=1}^{\min(\reff,\rhat)} \left( \theta_i^{2} + \rho_i^2  +  \dfrac{4 \rho_i }{ \theta_{i}^{2} D'_{\mu_X}(\rho_i)} \right)+ \sum_{i = \min(\reff,\rhat)+1}^{\max(r,\rhat)} (\theta_{i}^{2}+b^2),$$
where $\rho_{i} = D_{\mu_X}^{-1}(1/\theta_i)$ and we set $\theta_{i} =0$ for $i > r$.
Consequently, 
$$\mse{\weym{}} - \mse{\wopt{}} > \sum_{i = 1}^{\min(\reff,\rhat)}  \left(\rho_i + \dfrac{2}{\theta_i^2 D_{\mu_X}'(\rho_i)}\right)^2 + \sum_{i = \min(\reff,\rhat)+1}^{\max(r,\rhat)} b^2 > 0,  $$
almost surely.
\end{cor}
Corollary \ref{th:perf reff} reveals that the optimal estimator can realize a significant improvement in performance relative to the EYM estimator whenever $b = O(1)$ and $\reff < r$. The corollary highlights the importance of reliably estimating $\reff$ instead of $r$. Now, consider the rank regularized optimization problem
\be \label{eq:gwopt}
\gwopt{}(\rhat) :=  \argmin_{w \in \mathbb{R}^{q}, ||w||_{\ell_0} = \rhat } \fronorm{\sum_{i=1}^{r} \theta_i u_i v_{i}^{H} - \sum_{i=1}^{q} w_{i} \widehat{u}_{i} \widehat{v}_{i}^{H}}\ee
We characterize the exact solution next.
\begin{Th}[Optimal rank regularized solution]\label{th:soft thresholding}
For arbitrary integer $1 \leq \widehat{r} \leq q$, the solution to (\ref{eq:gwopt}) is given by
$$\gwopt{}(\rhat) =   \rlargest{\rhat} \left[  \{\Re (\sum_{j=1}^{r} \theta_j \uhatu{i}{j} \vvhat{j}{i} )_{+} \}_{i=1}^{q} \right],$$
where for $x \in \mathbb{R}_{+}^{q}$,   $\rlargest{\rhat}(x)$   returns a $q \times 1$ vector whose $\rhat$ non-zero elements equal the $\rhat$ largest entries of $x$ while the remaining entries are identically zero.
\end{Th}

\noindent We state a conjecture on the delocalization of the bulk singular vectors and characterize the asymptotic limit of (\ref{eq:gwopt}) next.

\begin{conj}[Complete delocalization of bulk singular vectors] \label{conj:full delocalization}
Define $q  = \min(m,n)$.  Assume that $ D_{\mu_X}'(b^+)=-\infty$  and that for all $1 \leq i \leq q$, where $i$ depends on $n$
$$\max_{i} \left( \sigma_i(X) - \sigma_{i+1}(X) \right) \leq  O\left( \dfrac{\log n ~{\rm factors}}{n}\right),$$
with very high probability.
Then, we have that for large enough $n$ and every  $i_n> \reff$  and $j = 1, \ldots r$
$$\max _{i,j} |\uhatu{i}{j}|  \leq O\left( \dfrac{\log n ~{\rm factors} }{n^{1/2}}\right)  \qquad ~{\rm and }  \qquad \max _{i,j} |\vvhat{j}{i}| \leq O\left( \dfrac{\log m ~{\rm factors}}{m^{1/2}}\right),$$
with high enough probability that we can establish their almost sure convergence to zero.
\end{conj}

\begin{cor}[Limiting rank regularized weights]\label{th:full rank soln}
Assuming Conjectures \ref{conj:partial delocalization} and \ref{conj:full delocalization} hold, we have that
$$\gwopt{i}(\rhat)  \convas
\begin{cases}
-2 \dfrac{D_{\mu_X}(\rho_i)}{D'_{\mu_X}(\rho_i)}  & \textrm{ for } i = 1, \ldots, \min(\reff,\rhat) \\ \\
0 & {\rm otherwise.}
\end{cases}$$
Consequently, even though, for finite $n$
$$\mse{\gwopt{}(q)} \leq \mse{\gwopt{}(\reff)} \leq \mse{\wopt{}(\reff)},$$
as $n \to \infty$ we have that
$$\mse{\gwopt{}(q)} - \mse{\gwopt{}(\reff)} \convas 0 \qquad {\rm and } \qquad  \mse{\gwopt{}(q)} - \mse{\wopt{}(\reff)} \convas 0.$$
\end{cor}

Corollary \ref{th:full rank soln} shows that when there is delocalization in the singular vectors then, in the large matrix limit, optimal performance is attained by estimating the effective rank $\reff$, applying shrinkage to the informative $\reff$ components and thresholding (to zero) the remaining components. In other words, there are vanishing (with $n$) performance losses when the coefficients given by $\wopt{}(\reff)$ are used in place of $\gwopt{}(q)$.
We believe that Conjectures \ref{conj:partial delocalization} and \ref{conj:full delocalization} hold in the signal-plus-noise matrix with missing entries setting considered in Section \ref{sec:missing data} so that Corollaries \ref{th:perf reff} and \ref{th:full rank soln} will apply there as well. This is pertinent because we now describe an algorithm for consistently estimating $\wopt{}$ \textit{directly from data} by exploiting the information in the singular value spectrum of the signal-plus-noise matrix.

\subsection{A new algorithm for improved denoising}\label{sec:algorithm}
\begin{algorithm}[t]
\centering
\caption{OptShrink: A new algorithm for low-rank matrix denoising by optimal, data-driven singular value shrinkage.}
\label{alg:new}
\begin{algorithmic}[1]\label{alg}
\STATE Input: $\wtX=  n \times m$ signal-plus-noise matrix
\STATE Input: $\rhat= \mbox{Estimate of the effective rank of the latent low-rank signal matrix}$
\STATE Compute $\wtX = \sum_{i=1}^{q} \widehat{\sigma}_{i} \uhat{i} \vhat{i}^H$
\STATE Compute $\widehat{\Sigma}_{\rhat} = \diag(\widehat{\sigma}_{\rhat+1}, \ldots \widehat{\sigma}_q) \in \mathbb{R}^{(n -\rhat) \times (m-\rhat)}$
\FOR {$i = 1, \ldots \rhat$}
\STATE  Compute $\widehat{D}(\widehat{\sigma}_i;\widehat{\Sigma}_{\rhat})$ using (\ref{eq:Dzhat})  and $\widehat{D}'(\widehat{\sigma}_{i};\widehat{\Sigma}_{\rhat})$ using (\ref{eq:Dpzhat}) \vspace{0.1cm}
\STATE{ Compute $\widehat{w}^{\rm opt}_{i,\rhat} = -2 \,\dfrac{\widehat{D}(\widehat{\sigma}_i;\widehat{\Sigma}_{\rhat})}{\widehat{D}'(\widehat{\sigma}_i;\widehat{\Sigma}_{\rhat})}$} \vspace{0.1cm}
\ENDFOR
\RETURN $\widehat{S}_{\rm opt} = \sum_{i=1}^{\rhat} \widehat{w}^{\rm opt}_{i,\rhat}\, \uhat{i} \,\vhat{i}^{H}$ = denoised estimate of the rank $\rhat$  signal matrix
\RETURN  (optional) Compute estimate of MSE using (\ref{eq:mse hat})
\RETURN  (optional) Compute estimate of relative MSE using (\ref{eq:rel mse hat})
\end{algorithmic}
\end{algorithm}

Equation (\ref{eq:wopt thm}) shows that the optimal estimator in the large matrix limit is given by
$$\wopt{i} = - 2 \dfrac{D_{\mu_X}(\rho_i)}{D'_{\mu_X}(\rho_i)}  + o(1),$$
where $\rho_i$ is the large matrix limit of the $i$-th largest singular value. In the finite $n,m$ setting, for $i = 1, \ldots, \reff$, $\widehat{\rho}_i= \widehat{\sigma}_i$ is a biased, but asymptotically consistent estimator of $\rho_i$. We now describe an algorithm for estimating $\wopt{i}$ using a single signal-plus-noise matrix.

For a matrix $X \in \mathbb{K}^{n \times m}$. Define
\begin{subequations}
 \be \label{eq:Dzhat}
\widehat{D}(z;X)  := \dfrac{1}{n} \Tr \left(z\, (z^2\,I - XX^H)^{-1} \right) \cdot \dfrac{1}{m} \Tr \left( z \,(z^2\,I - X^HX)^{-1}\right)
\ee
\textrm{and}
\begin{multline}\label{eq:Dpzhat}
\widehat{D}'(z;X) :=  \dfrac{1}{n} \Tr \left[z\, (z^2\,I - XX^H)^{-1} \right] \cdot \dfrac{1}{m} \Tr \left[ -2z^2 \,(z^2\,I - X^HX )^{-2} + \,(z^2\,I - X^HX )^{-1}\right]  \\
+ \dfrac{1}{m} \Tr \left[z\, (z^2\,I - X^H X)^{-1} \right] \cdot \dfrac{1}{n} \Tr \left[ -2z^2 \,(z^2\,I - XX^H )^{-2} + \,(z^2\,I - XX^H )^{-1}\right].
\end{multline}
\end{subequations}
By construction (and the definition of the $D$-transform), $\widehat{D}(z;X)  \convas D_{\mu_X}(z)$ and $\widehat{D}'(z;X) \convas D'_{\mu_X}(z)$ for $z$ outside the support of $\mu_X$. We now show how the spectrum of $\wtX$  can be used to estimate $\mu_X$. To that end, we establish a useful identify by first defining
\[
\mu_{X,\rhat} = \frac{1}{n-\rhat} \sum_{i=\rhat+1}^{n} \delta_{\si_{i}(X_{n})}.
\]
Then, it is easy to see that for fixed (with $n$) $\rhat$, $\mu_{X_{\rhat}} \convas \mu_{X,0}$.
Thus, if
$$\widehat{\Sigma}_{\rhat} = \diag(\widehat{\sigma}_{\rhat+1}, \ldots \widehat{\sigma}_q) \in \mathbb{R}^{(n -\rhat) \times (m-\rhat)}$$
is a diagonal matrix containing the $q - \rhat$ ``noise'' singular values of $\wtX$, then, by construction, and whenever $\widehat{\sigma}_{i} \convas \rho_i > b$, then $\widehat{D}(\widehat{\sigma}_i;\widehat{\Sigma}_{\rhat}) \convas D_{\mu_X}(\rho_i)$ and  $\widehat{D}'(\widehat{\sigma}_i;\widehat{\Sigma}_{\rhat}) \convas D'_{\mu_X}(\rho_i).$  Hence, we form a consistent estimate of $\wopt{i}$ as described in Algorithm \ref{alg:new}. The methods described in Section \ref{sec:rank est} can be used to form an estimate of $\rhat$.

By Theorem \ref{th:limiting mse}, we can compute an estimate of the absolute and relative mean squared error (defined as $\textrm{MSE}/|| S||_{F}^{2}$) as
\begin{subequations}
\be \label{eq:mse hat}
\widehat{\textrm{MSE}}_{\rhat} = {\sum_{i=1}^{\rhat}  \dfrac{1}{\widehat{D}(\widehat{\sigma}_i;\widehat{\Sigma}_{\rhat})}}  - {\sum_{i=1}^{\rhat} (\widehat{w}^{\rm opt}_{i,\rhat})^2}
\ee
\be \label{eq:rel mse hat}
\textrm{rel}\widehat{\textrm{MSE}}_{\rhat} = 1 - \dfrac{\sum_{i=1}^{\rhat} (\widehat{w}^{\rm opt}_{i,\rhat})^2}{\sum_{i=1}^{\rhat}  \dfrac{1}{\widehat{D}(\widehat{\sigma}_i;\widehat{\Sigma}_{\rhat})}},
\ee
\end{subequations}
respectively. A value for $\textrm{rel}\widehat{\textrm{MSE}}_{\rhat}$ near $0$ indicates very good low-rank signal matrix approximation while a value near $1$ indicates a poor approximation. These metrics might be better proxies for the noisiness of a signal-plus-noise matrix than the condition number or the spectral gap. We conclude with a statement of the theoretical consistency of the $\wopt{i}$ produced by Algorithm \ref{alg:new}.

\begin{Th}\label{th:wopt consistent}
Assume that $\rhat = \reff$. Then for $1 \leq i \leq \rhat$, we have that
$$\widehat{w}^{\rm opt}_{i,\rhat}  \convas -2 \dfrac{D_{\mu_X}(\rho_i)}{D'_{\mu_X}(\rho_i)} $$
\end{Th}
\begin{proof}
This is a straightforward consequence of Theorem \ref{th:wyem vs wopt}-a) and the fact that the almost sure limit of (\ref{eq:Dzhat}) leads (as described in the introduction of \cite{benaych2011svd}) directly to the $D$-transform.
\end{proof}

\section{Numerical Validation, Discussion and Extensions}\label{sec:discussion}

We now numerically validate our predictions. In the experiments that follow, we consider the model in (\ref{eq:model}) with $r = 1$, $n = m = 400$ and select $X$ to be an $n \times m$ matrix with i.i.d. $\mathcal{N}(0,1/m)$ entries.
For various values of $\theta$, Figure \ref{fig:comparison}-a) compares empirically computed $\wopt{1}$ averaged over $100$ trials with the (limiting) theoretical prediction given by the $p=1$ result in Theorem \ref{conj:missing}. Figure \ref{fig:comparison}-b) compares the realized normalized MSE  and shows that the EYM solution is near-optimal for large values of $\theta$ but far from sub-optimal for small values of $\theta$. The simulations validate the shrinkage-and-thresholding form of the solution for $\wopt{1}$ given by Theorem \ref{th:r1 soft thresholding} and show that Algorithm \ref{alg:new} realizes the predicted performance gains.

We now consider the optimization problem in (\ref{eq:wopt r}) and evaluate the performance of the various algorithms for various values of $\rhat$ for $\theta = 10$ and $\theta = 2$.  Here, $\reff = 1$ and Corollary \ref{th:perf reff} predicts that the optimal (oracle) algorithm should significantly outperform the EYM algorithm whenever $\rhat > \reff$. Figure \ref{fig:comparsion rank} shows the validity of this prediction and also shows that even though Algorithm \ref{alg:new} is suboptimal, relative to the oracle estimator, it is able to largely mitigate the effect of $\reff$ overestimation due to the shrinkage effect.

Figure \ref{fig:risk estimate} compares the normalized MSE estimates as a function of $\theta$, produced by Algorithm \ref{alg:new} to the empirical values for the setting where $\rhat = r = 1$ and $\theta_1 = \theta $ and where $\rhat = r = 2$, $\theta_1 = 20$ and $\theta_2 = \theta$.  As expected the estimates, produced are accurate whenever $\reff = \rhat$.

We now validate Theorem \ref{conj:missing}. We fix $r = 1$ and $\theta_1 = \theta = 2$ in (\ref{eq:missing model}) and vary $p$, the proportion of entries with missing data. We sample $u_1$ and $v_1$ uniformly at random from the unit hypersphere so that the low-coherence conditions in Theorem \ref{conj:missing} are met. Theorem \ref{conj:missing} predicts  that $\wopt{1} \to 0$ (asymptotically) when $p < \sqrt{n/m}/\theta^2 =  0.25$. Figure \ref{fig:comparison p} shows the accuracy of the prediction and the significant improvement in performance of the oracle estimator and Algorithm \ref{alg:new} relative to the EYM estimator.

\subsection{Suboptimality of singular value thresholding}

We now 	compare our algorithm to regularized matrix estimates obtained as the solution to the optimization problem
\be \label{eq:svt problem}
 \widehat{S}_{\rm svt, \lambda} = \argmin_{S} \fronormsq{\wtX - S} + 2 \lambda || S||_*,
\ee
where $|| \cdot ||_{*}$ is the nuclear norm (or the sum of the singular values of the argument matrix). The optimization problem in (\ref{eq:svt problem}) yields the closed-form solution \cite{cai2010singular}
\be  \label{eq:svt}
\widehat{S}_{\rm svt, \lambda}  = \Uhat \diag((\widehat{\sigma} - \lambda)_+) \Vhat^H.
\ee
The resulting singular value thresholded (SVT) matrix corresponds to the weighting
$$ w_{{\rm svt},i}(\lambda) = \begin{cases} \widehat{\sigma}_{i} - \lambda & {\rm if  }~ \widehat{\sigma}_{i} > \lambda \ \\
0 & {\rm otherwise}.
\end{cases}
$$

Figure \ref{fig:comparison shrink}-a) and b) compare the resulting soft-thresholding operator associated with the SVT approximation with the optimal and the EYM solutions for $\lambda = 1, 2$ as a function of $\theta$ and $\weym{}$, respectively for the same $r = 1$, $n = m$ setting in (\ref{eq:model}) with $X_{ij}$ i.i.d. $\mathcal{N}(0,1/m)$. Here $b = (1+\sqrt{c}) = 2$.

While SVT with $\lambda = 2$ can yield comparable shrinkage (in the small $\theta$ regime) and thresholding (below $\theta = 1$) as the optimal estimator, $w_{{\rm svt}}(2) - \wopt{}$ will be large for moderate $\theta$ so that by Theorem \ref{th:limiting mse}-c) we expect SVT to be  suboptimal for larger values of $\theta$.  Figure \ref{fig:svt versus opt} compares the performance of Algorithm 1 and the optimal estimator to the SVT algorithm with $\lambda = 1$ and $\lambda =2$. SVT is significantly suboptimal as expected. Our results show that our algorithm would outperform SVT with convex shrinkage functions for any of the general family of noise models considered here.

\subsection{Better singular value shrinkage with non-convex potential functions?}

A closer examination of Figure \ref{fig:comparison shrink}-a) and b) reveals that the optimal estimator shrinks less for larger values of $\theta$ than the SVT possibly can. In fact, the optimal estimator will generically yield a non-convex shrinkage function which scales as $$\wopt{i} \approx \widehat{\sigma}_i \left(1- O\left(\dfrac{1}{\widehat{\sigma}_i^2}\right)\right),$$ for large $\widehat{\sigma}_i$. Might singular value shrinkage with other non-convex potential functions generically outperform convex potential functions as well? These would be the non-convex analogs in the matrix setting of the non-negative Garrotte estimator \cite{breiman1995better} in the vector setting. Fully understanding their benefits and shortfalls, relative to Algorithm 1, remains an open line of inquiry.

\subsection{Role of informative components}
We conclude by reexamining the role of the principal (or leading) $\reff$ singular vectors of $\wtX$ in the solution of the optimization problem  (\ref{eq:gwopt}).  Theorem \ref{th:soft thresholding} shows that we should take the components $\uhat{i}$ and $\vhat{i}$ for which the inner product $\uhatu{i}{i}$ and $(\vhatv{i}{i})$ is $O(1)$. The supposition in (\ref{eq:wopt}) is that the principal components are these components.

However, in an expository paper by the author \cite{nadakuditi2013most}, it is shown that if the (limiting) spectrum of the noise-only matrix is supported on two disconnected intervals, then the middle components can be more informative than the principal components. Thus, while this work (via Theorem \ref{th:limiting mse}) brings into focus the importance of accurately estimating $\reff$, it is equally important to be able to identify the most informative components. The development of fast, accurate algorithms for the same for large matrix-valued datasets remains an important open problem.

\subsection{Extensions}

We have initial numerical evidence that the algorithm presented here outperforms the EYM estimator for the variety of applications described in \cite{chatterjee2012matrix}, even though they do not exactly fit the noise matrix models analyzed here. Extending the analysis of our algorithm to these models would shed further insight on the limits of low-rank signal matrix approximation.

We conclude by listing some directions of future research. These include 1) rigorously establishing the delocalization conjectures, 2) designing penalty functions that are robust to noise model mismatch, 3) clarifying the benefits, if any, of matrix regularization \cite{kakade2012regularization,deledalle2012risk,koltchinskii2011nuclear} with convex or non-convex penalty functions relative to rank regularized solutions for the unstructured low-rank signal matrix setting, 4) extending the methods developed to problems involving estimation of signal matrices with an unstructured low-rank  component and a sparse \cite{Candes:2011:RPC:1970392.1970395,chandrasekaran2009sparse,chandrasekaran2011rank,tao2011recovering,oymak2011finding} or  diagonal \cite{saunderson2012diagonal} component or low-rank structured component \cite{chu2003structured} and 5) developing minimax estimators, along the lines of the work in \cite{chatterjee2012matrix}, except for the more general class of noise models considered here.

Lastly, consider Theorem \ref{th:r1 soft thresholding}, where it is shown that for $\theta < 1/D_{\mu_X}(b^+)$, $\mse{\wopt{1}} \convas \theta^2$. In this regime, is there another (non-SVD based) algorithm that can estimate the signal matrix with mean-squared-error $\theta^2 -O(1)$?  More generally, is there a non-SVD based algorithm that can (reliably) recover the (unstructured) low-rank signal matrix in the regime where the SVD based methods break down? This is a largely open question whose answer would better clarify the interplay between the limits of SVD-based estimation of the signal matrix singular vectors and the fundamental limits of estimation of the signal matrix itself.  We leave these questions for future work.

\begin{figure}[t]
\subfloat[$\wopt{1}$ and $\weym{1}$ versus $\theta$.]{
\includegraphics[trim = 75 217 75 234, clip = true,width=5.05in]{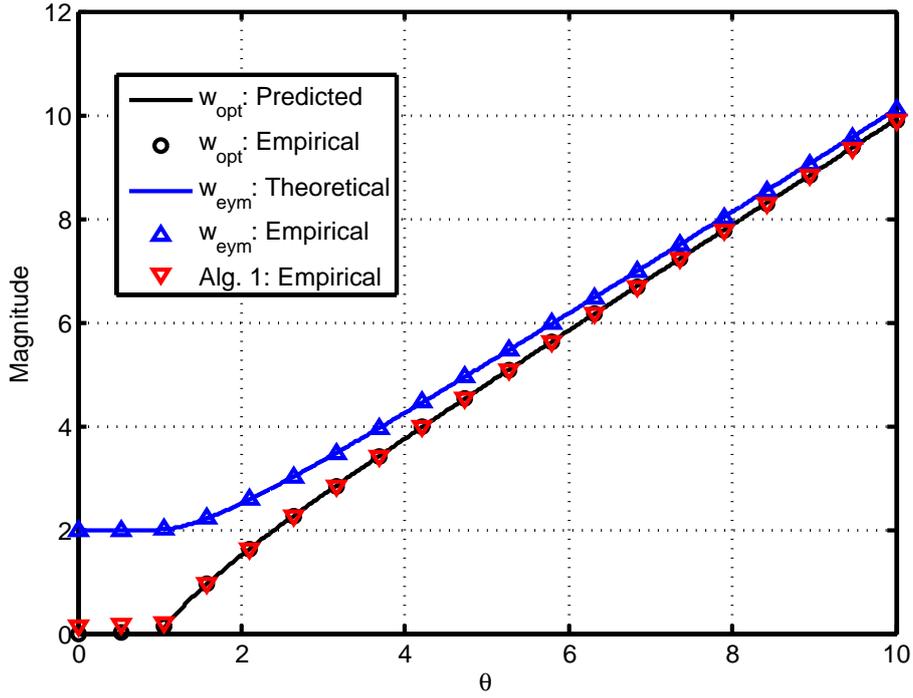}
}
\\
\subfloat[$\fronormsq{\theta u v^{H} - w \uhat{1} \vhat{1}^H}/\theta^2$ versus $\theta$.]{
\includegraphics[trim = 75 217 75 234, clip = true,width=5.05in]{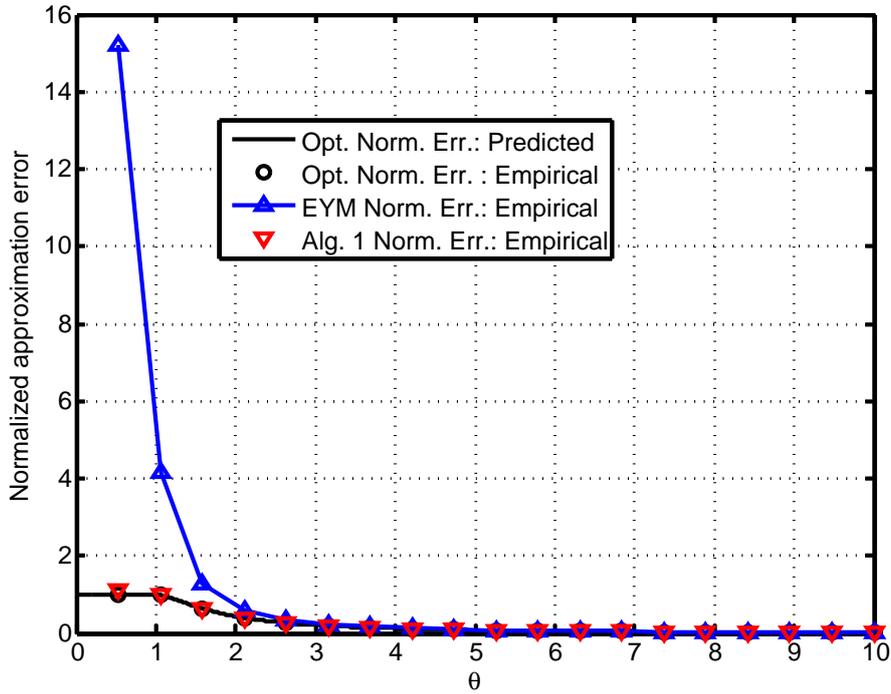}
}
\caption{For the model in (\ref{eq:model}) with $r = 1$, $n = m = 400$ and $X$ an $n \times m$ matrix with i.i.d. $\mathcal{N}(0,1/m)$ entries,  for various values of $\theta:=\theta_1$ , (a) we compare the theoretically predicted $\wopt{1}$  using (\ref{eq:wopt thm}) with the $\weym{1}$ computed  using (\ref{eq:weym thm}) (so that they precisely correspond to the $p=1$ prediction in Theorem \ref{conj:missing}) with empirically computed values of the same (averaged over $100$ trials).  Here we set $\rhat = 1$ in Algorithm \ref{alg:new}. (b) plots the realized normalized approximation errors and compares them to the (oracle) performance of the optimal detector predicted in Theorem \ref{th:limiting mse}. }\vskip-0.05cm
\label{fig:comparison}
\end{figure}

\begin{figure}[t]
\subfloat[Normalized MSE versus $\rhat$: $\theta = 10$.]{
\includegraphics[trim = 75 235 75 245, clip = true,width=5.25in]{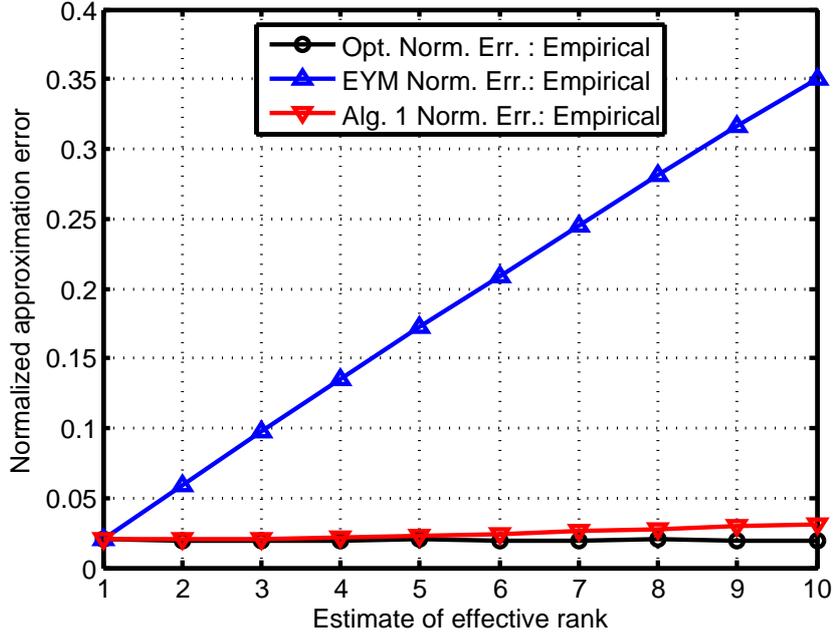}
}
\\
\subfloat[Normalized MSE versus $\rhat$: $\theta=2$.]{
\includegraphics[trim = 75 235 75 245, clip = true,width=5.25in]{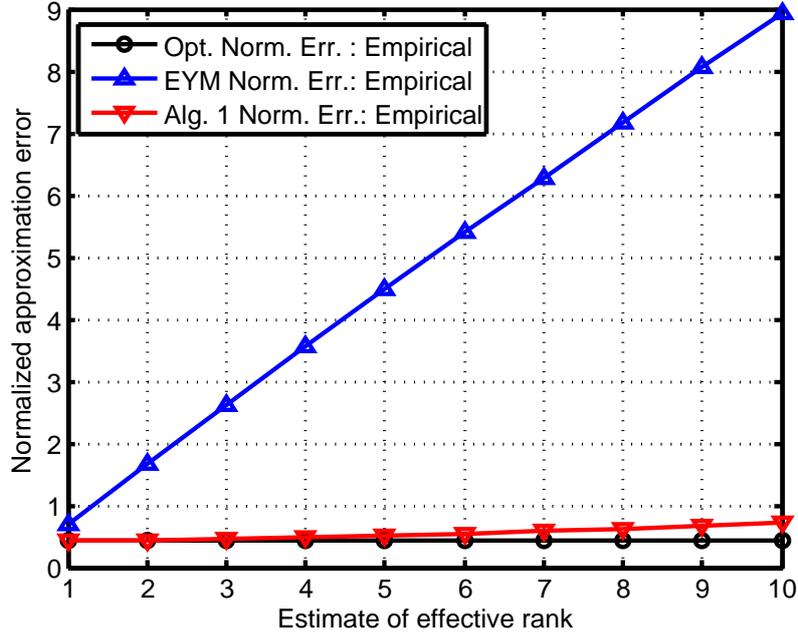}
}
\caption{Here, we are in the same setting as in Figure \ref{fig:comparison}, except, we evaluate the performance of Algorithm 1 and the EYM estimator of rank $\rhat$ to that of the rank $\rhat$ oracle optimal estimator computed using the left hand side of (\ref{eq:wopt thm}) for various values of $\rhat$. In a), $\theta = 10$ while in b) $\theta = 2$. }
\label{fig:comparsion rank}
\end{figure}

\begin{figure}
\centering
\subfloat[$\rhat = 1,$ $S= \theta u v^H$]{
\includegraphics[trim = 75 245 75 248, clip = true,width=5.25in]{{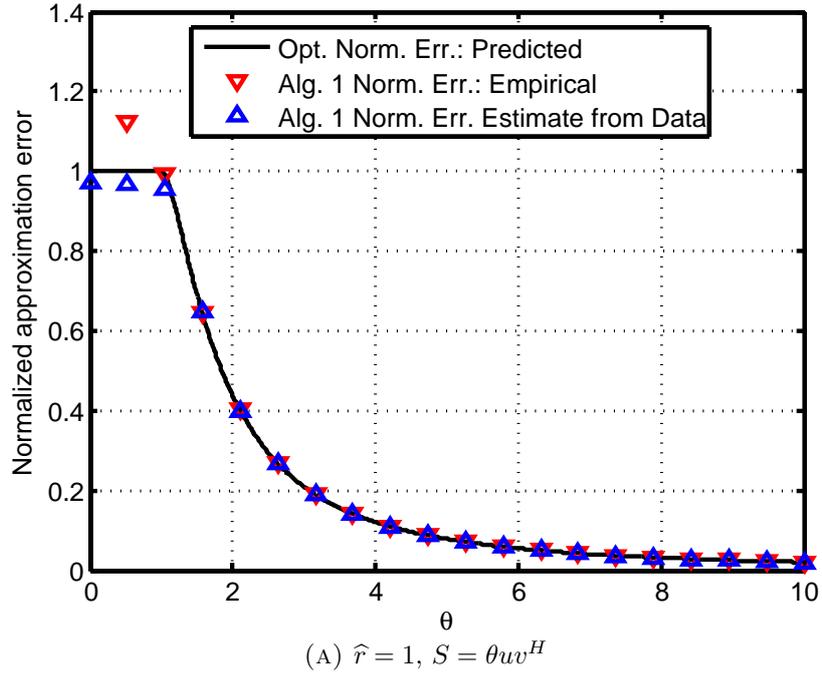}}}\\
\subfloat[$\rhat = 2$, $S= 20\, u_1 v_1^H + \theta u_2 v_2^H$]{
\includegraphics[trim = 75 245 75 245, clip = true,width=5.25in]{{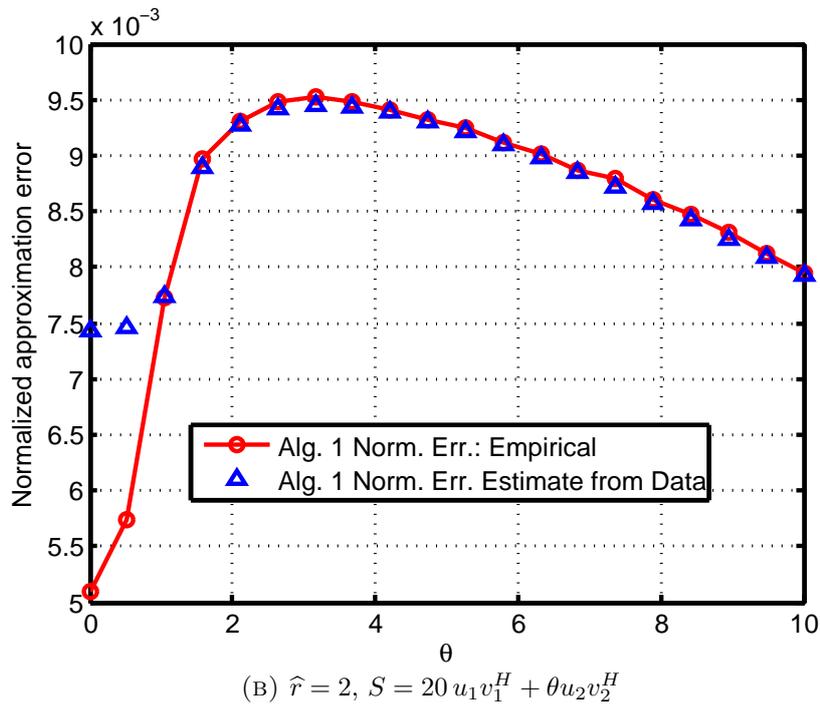}}
}
\caption{Here, we compare the normalized approximation error computed empirically with the estimate $ \textrm{rel}\widehat{\textrm{MSE}}_{\rhat}$ computed using (\ref{eq:rel mse hat}) as returned by Algorithm 1. When $\theta\leq 1$, $\reff = r-1$ so that one of the components becomes uninformative so that including it in the estimate will increase the realized error.}\vskip-0.05cm
\label{fig:risk estimate}
\end{figure}

\begin{figure}[t]
\subfloat[$\wopt{1}$ and $\weym{1}$ versus $p$.]{
\includegraphics[trim = 75 235 75 245, clip = true,width=5.25in]{{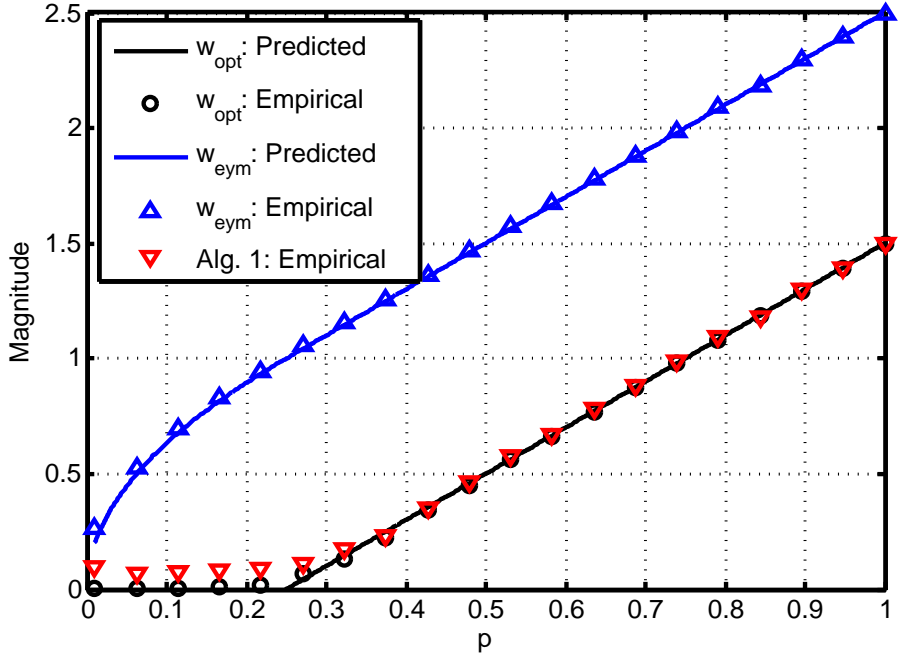}}
}
\\
\subfloat[$\fronormsq{p \theta u v^{H} - w \uhat{1} \vhat{1}^H}/p^2/\theta^2$ versus $p$; here $\theta = 2$.]{
\includegraphics[trim = 75 235 75 245, clip = true,width=5.25in]{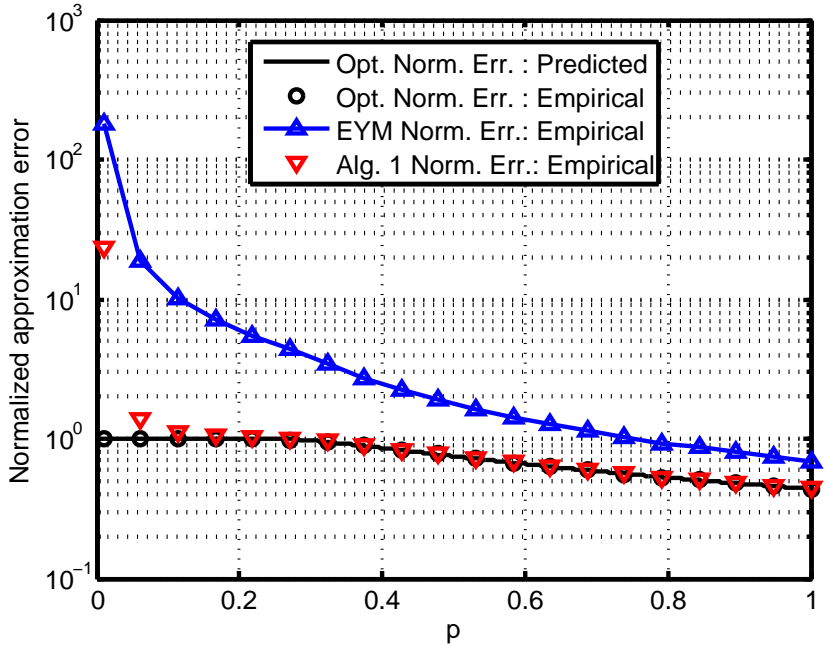}
}
\caption{For the model in (\ref{eq:missing model}), with $r = 1$, $n = m = 400$ and $X$ an $n \times m$ matrix with i.i.d. $\mathcal{N}(0,1/m)$ entries, we perform the same comparisons as in Figure \ref{fig:comparison} (averaged over $100$ trials), except we fix $\theta = 2$ and instead vary $p$, the proportion of entries with missing data. We sample $u_1$ and $v_1$ uniformly at random from the unit hypersphere so that the low-coherence conditions in Theorem  \ref{conj:missing} are met. Theorem  \ref{conj:missing} predicts that $\wopt{} \to 0$ (asymptotically) when $p < \sqrt{n/m}/\theta^2 =  0.25$.}
\label{fig:comparison p}
\end{figure}

\begin{figure}[t]
\subfloat[Shrinkage and thresholding operators as a function of $\theta$.]{
\includegraphics[trim = 75 235 75 245, clip = true,width=5.25in]{{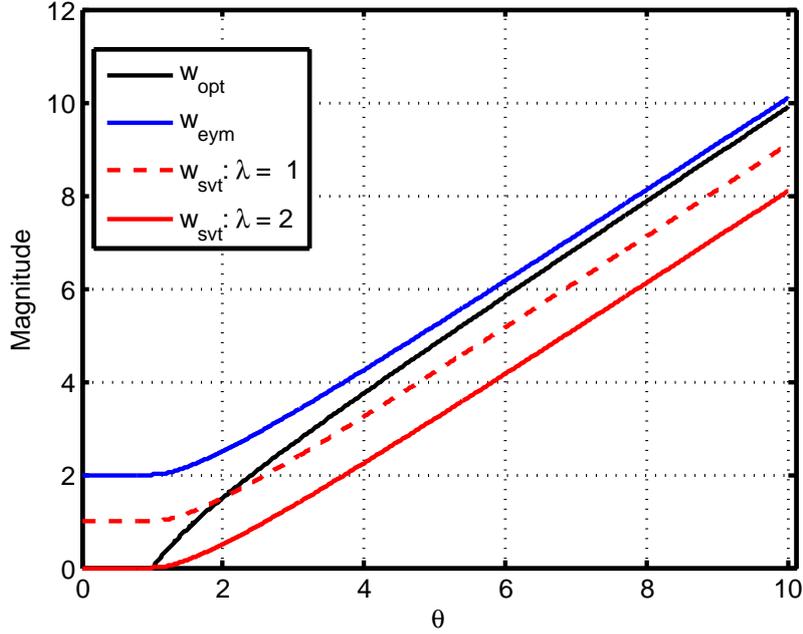}}
}
\\
\subfloat[Shrinkage and thresholding operators as a function of $\weym{}$.]{
\includegraphics[trim = 75 235 75 245, clip = true,width=5.25in]{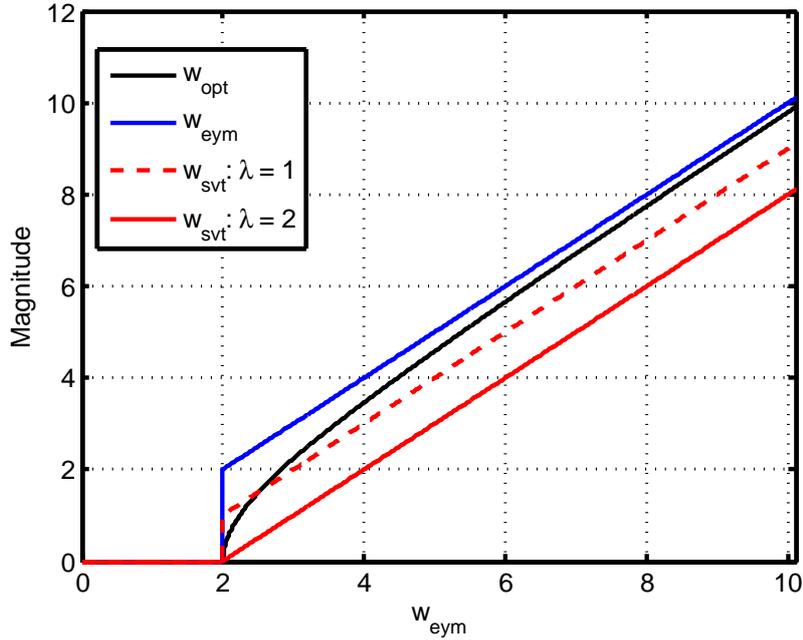}
}
\caption{Here we are in the same setting as Figure \ref{fig:comparison}. We plot $\wopt{}$, $\weym{}$ and $w_{{\rm svt},\lambda}$ for $\lambda = 1, 2$ as a function of $\theta$ and $\weym{}$. Note the non-convex nature of the shrinkage portion of the optimal shrinkage-and-thresholding operator, the optimality of the EYM solution for large values of $\theta$ (high SNR regime) and the sub-optimality of the SVT solution.}
\label{fig:comparison shrink}
\end{figure}

\begin{figure}
\centering
\includegraphics[trim = 75 235 75 245, clip = true,width=5.25in]{{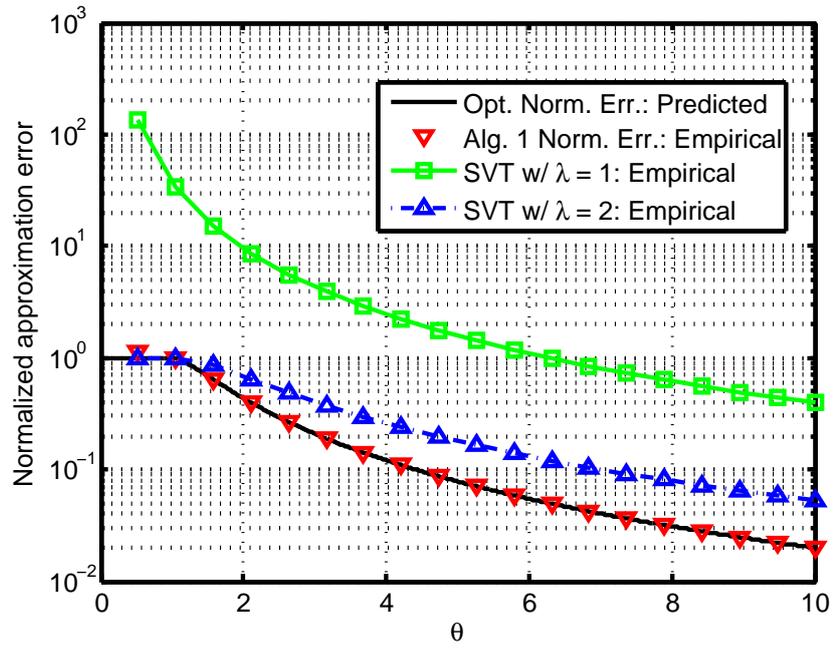}}
\caption{For the same setting as in Figure \ref{fig:comparison}, we compare the performance of Algorithm 1 (with $\rhat  =1$) with that of the SVT estimator in (\ref{eq:svt}) for $\lambda =1, 2$ .}\vskip-0.05cm
\label{fig:svt versus opt}
\end{figure}

\clearpage
\section{Proof of Theorems \ref{th:wyem vs wopt}, \ref{th:r1 soft thresholding} and \ref{th:soft thresholding}}

We first prove Theorem \ref{th:wyem vs wopt} -b).  Since $\weym{i} = \widehat{\sigma}_{i}$, Theorem \ref{th:wyem vs wopt}-b) follows immediately from Theorem 2.9 in \cite{benaych2011svd}.  Next, we prove the first part of Theorem \ref{th:wyem vs wopt}-a) by showing that
$$\wopt{i} =  \left(\Re\{ \sum_{j=1}^{r} \theta_j \uhatu{i}{j} \vvhat{j}{i} \}\right)_{+}.$$ Theorem \ref{th:soft thresholding} follows by adopting the exact same approach, with some minor modifications so we shall omit its proof. We first establish some intermediate results.

\begin{lemma}\label{lemma:diag approx}
Let $A \in \mathbb{K}^{n \times m}$ and $q = \min(n,m)$. Consider the optimization problem
$$D_{opt} := \argmin_{D = \diag(\{d_1, \ldots, d_q\}), d_i \in \mathbb{R}_+} \fronorm{A - D},$$
where $\diag(\cdot)$ denotes a matrix with the arguments on the diagonal and zeros elsewhere (even for a rectangular matrix). Then
$$( D_{opt})_{ii} = \max(0,\Re( A_{ii} )).$$
\end{lemma}
\fl\begin{proof}
We first solve the unconstrained problem
$$D_{opt} := \argmin_{D = \diag(\{d_1, \ldots, d_q\})} \fronorm{A - D},$$
Note that
$$ \fronormsq{A - D} = \sum_{i=1}^{q} (A_{ii} - d_{i})^2 + \underbrace{\sum_{i \neq j} A_{ij}^2}_{\textrm{constant}} \geq \sum_{i \neq j} A_{ij}^2,$$
so that setting $d_{i} = A_{ii}$ attains the lower bound.  The additional constraint that $d_{i} \in \mathbb{R}_+$ yields the stated result which is simply a projection onto $\mathbb{R}_+$.
\end{proof}

\begin{corollary}\label{corr:diagonal}
For fixed $r$, the solution to the optimization problem
$$D_{opt} := \argmin_{D = \diag(\{d_1, \ldots d_r, 0, \ldots , 0 \}), d_i \in \mathbb{R}_+ } \fronorm{A - D},$$
is given by
$$ (D_{\rm opt})_{ii} = \max(0, A_{ii}) \qquad {\rm for } i = 1, \ldots r.$$
\end{corollary}
\fl Now consider the optimization problem
$$\wopt{} = \argmin_{w \in \mathbb{R}_{+}^{r}}  \fronorm{\sum_{i=1}^{r} \theta_i u_i v_i^{H} - \sum_{i=1}^{r} w_{i} \uhat{i} \vhat{i}^H }.$$
Let $U_{r} = \begin{bmatrix} u_{1} & \ldots & u_{r} \end{bmatrix}$, $V_{r} = \begin{bmatrix} v_{1} & \ldots & v_{r} \end{bmatrix}$, $\Theta_{r} = \diag(\theta_{1}, \ldots, \theta_{r})$, $\Uhat = \begin{bmatrix} \uhat{i} & \ldots \uhat{n} \end{bmatrix}$ and $\Vhat = \begin{bmatrix} \vhat{i} & \ldots \vhat{m} \end{bmatrix}$. Then for $W = \diag(w_{1}, \ldots, w_{r}, 0, \ldots 0)$, the optimization problem can be rewritten as

$$\wopt{}=  \argmin_{w \in \mathbb{R}_{+}^{r}, W = \diag(w)}  \fronorm{ U_{r} \Theta_{r} V_{r}^H -  \Uhat W \Vhat^{H}}.$$
By the unitary invariance of the Frobenius norm we have that
$$\fronorm{ U_{r} \Theta_{r} V_{r}^{H} - \Uhat W \Vhat^{H}} = \fronorm{\Uhat^{H}U_{r} \Theta_{r} V_{r}^{H} \Vhat -  W }.$$
Let $K =\Uhat^{H}U_{r} \Theta_{r} V_{r}^{H} \Vhat$. Then,

\begin{align*}
K &= \begin{bmatrix} \uhat{1}^{H}u_1  &  \ldots & \uhat{1}^H u_{r}\\
                                \vdots                       &   \vdots                   & \vdots                                \\
                                \uhat{r}^H u_1 & \ldots & \uhat{r}^H u_{r}  \\ \\
                                \uhat{r+1}^H u_1 &  \ldots & \uhat{r+1}^H u_{r}  \\
                                \vdots                       &   \vdots                   & \vdots                                \\
                                \uhat{n}^H u_1 & \ldots & \uhat{n}^H u_{r}  \\
\end{bmatrix}
\begin{bmatrix}
\theta_{1} &  &  \\
                 & \ddots  &  \\
                &               & \theta_{r}
\end{bmatrix}
\begin{bmatrix} v_{1}^H \vhat{1} & \ldots & v_{1}^H \vhat{r} & v_{1}^H \vhat{r+1} & \ldots & v_{1}^H \vhat{m} \\
                              \ldots                       &   \ldots                   &     \ldots                          & \ldots    	  \\
  v_{r}^H \vhat{1} & \ldots & v_{r}^H \vhat{r} &  v_{r}^H \vhat{r+1} & \ldots & v_{r}^H \vhat{m}\\
  \end{bmatrix} \\[0.25cm]
& =
 \sum_{j= 1}^{r}  \theta_j \begin{bmatrix} \uhat{1}^{H}u_j \\
                                \vdots                    \\
                                \uhat{r}^H u_j \\ \\
                                \uhat{r+1}^H u_j   \\
                                \vdots                       \\
                                \uhat{n}^H u_j  \\
\end{bmatrix}
\begin{bmatrix} v_{j}^H \vhat{1} & \ldots & v_{j}^H \vhat{r} & v_{j}^H \vhat{r+1} & \ldots & v_{j}^H \vhat{m} \\ \end{bmatrix} \\
\end{align*}
Expanding out the diagonal entries of $K$ we get
\begin{align*}
K & = \sum_{j=1}^{r}
\begin{bmatrix}
\theta_{j} \left(\uhat{1}^{H} u_j \right) \cdot \left(v_{j}^{H} \vhat{1}\right)  &  *       &* \\
					*						        & \ddots  & *   \\
                                                     *                                                                       &  *           & \theta_{j} \left(\uhat{m}^{H} u_j \right) \cdot \left(v_{j}^{H} \vhat{n}\right)\\
\end{bmatrix}\\
\end{align*}
so that (deterministically),
\be \label{eq:Kii general}
K_{ii} =  \sum_{j=1}^{r} \theta_j \uhat{i}^{H} u_{j}  v_{j}^H \vhat{i}
\ee
and the solution
$$\wopt{i} = \max(0,\Re \sum_{j=1}^{r} \theta_j \uhat{i}^{H} u_{j}  v_{j}^H \vhat{i}) = (\Re \sum_{j=1}^{r} \theta_j \uhat{i}^{H} u_{j}  v_{j}^H \vhat{i})_{+}  , $$
follows immediately from (\ref{eq:Kii general}) by the application of Corollary \ref{corr:diagonal}. We have thus proved the equality on the left-hand side of  Theorem \ref{th:wyem vs wopt}-a). It is easy to see how this approach yields Theorem \ref{th:soft thresholding}.

We now prove the limit characterization portion of Theorem \ref{th:wyem vs wopt}-a). In  \cite[Theorem 2.10 c)]{benaych2011svd}, it was proved that for $j =1, \ldots, r,$  and $i \neq j$ such that $\theta_{i}^2 > 1/D_{\mu_X}(b^+)$, $\uhat{i}^{H} u_{j} \convas 0$ and  $ v_{j}^H \vhat{i} \convas 0$. Consequently,
$$K_{ii} =  \sum_{j=1}^{r} \theta_j\, \uhat{i}^{H} u_{j}  \cdot v_{j}^H \vhat{i} \convas \theta_{i}\, \uhat{i}^{H} u_{i}  v_{i}^H \vhat{i}.$$
Let $\rho_{i} = D^{-1}_{\mu_X}(1/\theta_i)$. In  \cite[Theorem 2.10 c)]{benaych2011svd} it was shown that
$$ |\uhat{i}^{H} u_{i}|^{2} \convas \dfrac{-2 \phi_{\mu_X}(\rho_i)}{\theta_{i}^2 D'_{\mu_X}(\rho_i)} \qquad \textrm { and } \qquad  |\vhat{i}^{H} v_{i}|^{2} \convas \dfrac{-2 \phi_{\widetilde{\mu}_X}(\rho_i)}{\theta_{i}^2 D'_{{\mu}_X}(\rho_i)},$$
where $\widetilde{\mu}_X=c\mu_X+(1-c)\delta_0$ and for any \pro measure $\mu$,
\be\label{9710.19h07}\phi_\mu(z):=\int \f{z}{z^2-t^2}\ud \mu(t).\ee
While there is ambiguity in the sign (or phase, when complex valued) of the individual singular vectors, the proof in \cite{benaych2011svd} shows that
\be \label{eq:pip analytical}
\uhat{i}^{H} u_{i}  \, v_i^H \vhat{i}  \convas \sqrt{ \dfrac{-2 \phi_{\mu_X}(\rho_i)}{\theta_{i}^2 D'_{\mu_X}(\rho_i)}\cdot \dfrac{-2 \phi_{\widetilde{\mu}_X}(\rho_i)}{\theta_{i}^2 D'_{{\mu}_X}(\rho_i)}} = \dfrac{2 \sqrt{ \phi_{\mu_X}(\rho_i) \cdot  \phi_{\widetilde{\mu}_X}(\rho_i) }}{ \theta_{i}^{2} D'_{\mu_X}(\rho_i)} .
\ee
However, $D_{\mu_X}(z)   = \phi_{\mu}(z) \cdot \phi_{\widetilde{\mu}}(z)$ so that $\phi_{\mu}(\rho_i) \cdot \phi_{\widetilde{\mu}}(\rho_i) = D_{\mu_X}(\rho_i)  = D_{\mu_X}(D^{-1}_{\mu_X}(1/\theta_{i}^{2})) = 1/\theta_i^2$, so that
\begin{equation}\label{eq:key limit}
  \theta_{i} (\uhat{i}^{H} u_{i})   ( v_i^H \vhat{i})  \convas  \dfrac{-2 }{ \theta_{i}^{2} D'_{\mu_X}(\rho_i)} = -2 \dfrac{D_{\mu_X}(\rho_i)}{D'_{\mu_X}(\rho_i)}.
  \end{equation}
This gives the limit on the right hand side of part a).

To prove part c), we note that $\weym{i} > \theta_{i}$ (as a consequence of Horn's interlacing inequalities \cite{hj91}) while, for large enough $n$, $\wopt{i} = \theta_i \uhat{i}^{H} u_{i} v_{i}^H \vhat{i} + o(1) < \theta_i$. Thus $\weym{i} > \wopt{i}$ for large enough $n$. Since $\uhat{i}^H u_{i} v_{i}^H \vhat{i}  \to 1$ for $\theta_i \to \infty$, $\weym{i} \convas \wopt{i}$ as $\theta_i \to \infty$.

We now prove Theorem \ref{th:r1 soft thresholding}. Note that when $r =1$,
$$\wopt{1} = \left(\Re \theta_1 \uhat{1}^{H} u_{1}  v_{1}^H \vhat{1}\right)_+.$$
When $r =1$ and $\theta_{1}^2 \leq 1/D_{\mu_X}(b^+)$ and $D'_{\mu_X}(b^+) = -\infty$, then by Theorem 2.11 of \cite{benaych2011svd}, $\uhat{1}^{H} u_{1} \convas 0$ and $v_{1}^{H}\vhat{1} \convas 0$. Consequently, $\wopt{1} \convas 0$ and we have established the phase transition (or shrinkage-and-thresholding form) of $\wopt{1}$ in Theorem \ref{th:r1 soft thresholding}. The expressions for $\mse{\wopt{}}$ and $\mse{\weym{}}$ are a straightforward consequence of Theorem \ref{th:limiting mse}.

\section{Proof of Theorems \ref{th:limiting mse} and Corollaries \ref{th:perf reff} and \ref{th:full rank soln}}
Here, we have that
\begin{align*}
\mse{w} &=  \fronormsq{\sum_{i=1}^{r} \theta_{i} u_i v_i^H - \sum_{i=1}^{r} w_{i} \uhat{i}\vhat{i}^{H}}\\
& = \sum_{i}  \theta_i^2 + \sum_{j} w_j^2 - 2 \Re \,\textrm{Tr} \sum_{i,j} \theta_i w_j u_i v_i^{H} \vhat{j}\uhat{j}^{H} \\
& = \sum_{i} \theta_i^2 + \sum_i w_i^2 - 2\, \textrm{Tr} \sum_{i} \theta_i w_i  u_i v_i^{H} \vhat{i} \uhat{i}^H   - 2 \, \Re \, \textrm{Tr} \sum_{i \neq j} \theta_i w_j u_i v_i^{H} \vhat{j} \uhat{j}^{H} \\
\end{align*}
In  \cite[Theorem 2.10 c)]{benaych2011svd}, it was proved that for $j =1, \ldots, r,$  and $i \neq j$ such that $\theta_{i}^2 > 1/D_{\mu_X}(b^+)$, $\uhat{i}^{H} u_{j} \convas 0$ and  $ v_{j}^H \vhat{i} \convas 0$. Hence,
\begin{align*}
\mse{w} &= \sum_i (\theta_i^2 + w_i^2 - 2  \theta_i w_i  \uhat{i}^{H} u_{i} \cdot  \vhat{i}^{H} v_{i} ) - 2 \, \Re \sum_{i \neq j}  \theta_i w_j  \underbrace{\uhat{j}^H u_i}_{\convas 0}  \underbrace{v_i^{H} \vhat{j}}_{\convas 0} \\
& \convas \sum_{i=1}^{r} \left( \theta_i^2 + \dfrac{4 w_i }{ \theta_{i}^{2} D'_{\mu_X}(\rho_i)} + w_i^2 \right), \\
\end{align*}
where we have substituted (\ref{eq:key limit}) to give us the final expression in the stated result.

Theorem \ref{th:limiting mse}-a) and b) follow from substituting the limiting values of $\wopt{i}$ and $\weym{i}$ given by Theorem \ref{th:wyem vs wopt} in the derived expression. Theorem \ref{th:limiting mse}-c) follows easily by simple algebraic manipulation of the limiting expressions for $\mse{w}$ and $\mse{\wopt{}}$.  The portions of Corollaries \ref{th:perf reff} and \ref{th:full rank soln} that characterize the structure of the limiting weights follows immediately from Conjecture \ref{conj:partial delocalization} and Conjecture \ref{conj:full delocalization} via an application of Theorem \ref{th:soft thresholding}.

We now consider the asymptotic squared error. Note that
$$  \sqrt{\mse{\wopt{}(\reff)}} - \sqrt{\mse{\gwopt{}}} =  \fronorm{ \sum_{i=1}^{r} \theta_{i} u_i v_i^H - \sum_{i=1}^{q} \gwopt{i}\uhat{i}\vhat{i}^{H}} - \fronorm{\sum_{i=1}^{r} \theta_{i} u_i v_i^H - \sum_{i=1}^{\reff} \wopt{i}\uhat{i}\vhat{i}^{H}}.
$$
By the triangle inequality we have
$$  \sqrt{\mse{\wopt{}(\reff)}} -\sqrt{ \mse{\gwopt{}}}
\leq \fronorm{\sum_{i =1}^{\reff}(\gwopt{i}-\wopt{i}) u_i v_i^H} + \fronormsq{\sum_{i=\reff+1}^{q} \gwopt{i}\uhat{i}\vhat{i}^{H}}.
$$
Since we have just shown that $\gwopt{i} \convas \wopt{i}$ for $i = 1, \ldots, \reff$, we have
$$ \fronormsq{\sum_{i =1}^{\reff}(\gwopt{i}-\wopt{i}) u_i v_i^H} \convas 0.$$
If we can show that
$$ \fronormsq{\sum_{i=\reff+1}^{q} \gwopt{i}\uhat{i}\vhat{i}^{H}} \convas 0,$$
then we can conclude that $\mse{\gwopt{}} - \mse{\wopt{}(\reff)} \convas 0$ and we are done. To that end, we shall utilize the claim from Conjecture \ref{conj:partial delocalization} that the $o(n)$ leading coefficients of $\gwopt{i}$  corresponding to the edge (or principal) singular vectors will be bounded by $O(\log n ~{\rm factors}/n^{1/3})$  and the claim from Conjecture \ref{conj:full delocalization} that $O(n)$ of $\gwopt{i}$ coefficients corresponding to the bulk singular vectors will be bounded by $O(\log n ~{\rm factors}/n)$ with very high probability. This gives us
\begin{align*}
\fronormsq{  \sum_{i=\reff+1}^{q}  \gwopt{i}\uhat{i}\vhat{i}^{H}} & =
 \sum_{i > \reff, i \in {\rm bulk}} (\gwopt{i})^2 + \sum_{i> \reff, i \in {\rm edge}}  (\gwopt{i})^2  \\
 & \leq  O(n) O\left( \dfrac{\log n ~{\rm factors}}{n^2}\right)+ o(n) O\left( \dfrac{\log n ~{\rm factors}}{n^{2/3}}\right) \\
 & \leq \underbrace{O\left( \dfrac{\log n ~{\rm factors}}{n}\right) + O\left( \dfrac{\log n ~{\rm factors}}{n^{2/3}}\right)}_{\convas 0}.
 \end{align*}
If the probability is high enough we will be able to conclude that $\mse{\gwopt{}(\reff)} \convas \mse{\wopt{}(\reff)}$. Repeating this calculation with $\wopt{}(\rhat)$ and utilizing Conjecture \ref{conj:partial delocalization}  gives us the expression for the asymptotic squared error in Corollary \ref{th:perf reff}.

\section{Proof of Theorem \ref{conj:missing}}\label{proof:missing}
\noindent We begin by recalling that
\begin{equation}\label{eq:missing data}
\wtX = (\underbrace{U \Theta V^H}_{=:S} + X) \odot M = S \odot M + X \odot M,
\end{equation}
where $X$ is the noise-only matrix with $\E[X_{ij}] = 0$ and $\Var[X_{ij}]=1/m$ and
$$M_{ij} = \begin{cases} 1 & \textrm{ with probability } p \\
0 & \textrm{ with probability } 1- p.
\end{cases}
$$
Note that $\mathbb{E}_{M}[S \odot M] = p\,S$, so that (\ref{eq:missing data}) can be rewritten in a signal-plus-noise-plus-small-perturbation form\footnote{Thanks to Brendan Farrell for suggesting this approach to analyzing the problem.} given by
\begin{equation}\label{eq:missing equivalent}
 \wtX =\underbrace{\mathbb{E}[S \odot M]}_{= p\,S} +Z+\underbrace{\left(S \odot M - \mathbb{E}[S \odot M]\right)}_{=: \deltaS},
 \end{equation}
where ${Z}$ is the noise-only random matrix with missing entries given by
\begin{equation}\label{eq:alt construction}
{Z}_{ij} =
\begin{cases}
X_{ij} & \textrm{ with probability } p \\
0  & \textrm{ with probability } 1-p.
\end{cases}
\end{equation}
Let
\begin{equation}\label{eq:missing 1}
\overline{X} = p\,S + Z,
\end{equation}
so that, from (\ref{eq:missing equivalent}), $\widetilde{X} = \overline{X}+ \deltaS$. Let $\overline{X} = \sum_{i} \overline{\sigma}_{i} \ubar{i} \vbar{i}^{H}$ be the SVD of $\overline{X}$.  In lieu of (\ref{eq:wopt r missing}), consider the slightly modified optimization problem
\be \label{eq:wopt r missing alt}
{\woptbar{}} := \argmin_{w = [w_1 ~~ \cdots w_{r}]^T \in \mathbb{R}^{r}_{+}}  \fronormsq{\sum_{i=1}^{r} p \, \theta_i u_i v_i^H  - \sum_{i=1}^{r} w_{i} \ubar{i}\vbar{i}^{H}}.
\ee

We will first show that $\woptbar{i}$ is characterized by the stated expression in Theorem \ref{conj:missing}. Then we will show that $\sigma_{1}(\deltaS) \convas 0$, which we will utilize to prove that $\wopt{i} \convas \woptbar{i}$.

Comparing (\ref{eq:wopt}) to (\ref{eq:wopt r missing alt}) reveals that the left hand side of Theorem \ref{th:wyem vs wopt}-a) still holds except with $\theta_i \longmapsto p \,\theta_i$.  Consequently,
\begin{equation}\label{eq:woptbar}
\woptbar{i} =  \left(\Re\{ \sum_{j=1}^{r} p\, \theta_j \ubaru{i}{j} \vvbar{j}{i} \}\right)_{+}.
\end{equation}
We now establish the almost sure limit of the right hand side of (\ref{eq:woptbar}).

To that end, we first note that since $\mathbb{E}[X_{ij}] = 0$ and $\mathbb{E}[X_{ij}^{2}] = 1/m$, from (\ref{eq:alt construction}), we have that $\mathbb{E}[Z_{ij}] = 0$ and $\mathbb{E}[Z_{ij}^{2}] = p/m$. Moreover, since the higher order moments of the entries of $X$ were  assumed to be bounded, the higher order moments of the entries of $Z$ will be bounded as well. Consequently, it can be shown \cite{bai2010spectral} that
\begin{equation}\label{eq:mu z}
\ud \mu_{Z_n}(x) \convas \ud \mu_Z(x)  = \dfrac{\sqrt{4\,p^2\,c-(x^2-p-p\,c)^2}}{\pi \,p \, c\, x}\one_{(a,b)}(x) dx + \max \left(0,1-\dfrac{1}{c} \right)\delta_{0} ,
 \end{equation}
where $a=\sqrt{p}(1-\sqrt{c})$ and $b=\sqrt{p}(1+\sqrt{c})$ are the end points of the support of $\mu_Z$.  Here, $\mu_Z$ is the  famous Mar\v{c}enko-Pastur distribution  \cite{marchenko1967distribution}.  It is known \cite{bai2010spectral}, that $\sigma_{1}(Z) \convas b =\sqrt{p}(1+\sqrt{c})$.  Moreover, from the results of Bloemendal et al \cite[Theorems 2.4 and 2.5]{BEKYY13}, we have that for any $\{ u_i\}_{i=1}^{r}$ and $\{v_{i} \}_{i=1}^{r}$, independent of $Z$,
\begin{subequations}\label{eq:quad convergence}
\begin{equation}\label{eq:quad convergence u}
u_{i}^{H} (w^2 I_n - Z Z^{H})^{-1} u_{j} \convas \int \dfrac{\textrm{d}\mu_Z(t)}{w^2 - t^2} \, \delta_{ij}
\end{equation}
\text{and}
\begin{equation}\label{eq:quad convergence v}
v_{i}^{H} (w^2 I_m - Z^H Z)^{-1} v_{j} \convas \int \dfrac{\textrm{d}\mu_{\widetilde{Z}}(t)}{w^2 - t^2}\, \delta_{ij},
\end{equation}
\end{subequations}
where $\mu_{\widetilde{Z}} = c \mu_Z + (1-c) \delta_0$ (when $c<1$). An inspection of the proofs in \cite{benaych2011svd} reveals that the almost sure limits of these bilinear forms determine the almost sure limits of  $\sigma_{i}(\overline{X})$ and $\ubaru{i}{j}$ and $\vvbar{j}{i}$ for $i=1, \ldots, r$. Equation (\ref{eq:quad convergence}) asserts that these limits are the same as the limits that we would have obtained if $Z$ were i.i.d. Gaussian (and hence bi-unitarily invariant) with matching mean and variance as the $Z$ in (\ref{eq:alt construction}).  Consequently, the almost sure limit of $\woptbar{i}$ in (\ref{eq:woptbar}) will be the same as though $Z$ were i.i.d. Gaussian with mean zero and variance $p/m$  entries.  Hence,  by Theorem \ref{th:wyem vs wopt}-b)
$$\overline{w}_{i}^{\rm eym}:=\sigma_{i}(\overline{X}) \convas
\begin{cases}
\rho_i = D_{\mu_Z}^{-1}(1/p^{2} \theta_i^2) & {\rm if } ~p^2\,\theta_i^{2} > \dfrac{1}{D_{\mu_Z}(b^{+})} = p\, \sqrt{c}\\  \\
\sqrt{p}\,(1+\sqrt{c}) & {\rm otherwise},\\
\end{cases}
$$
while by Theorem \ref{th:wyem vs wopt}-a),
$$\woptbar{i} \convas
-2 \dfrac{D_{\mu_Z}(\rho_i)}{D'_{\mu_Z}(\rho_i)} \qquad  {\rm  if } ~ \theta^{2}_i  > \dfrac{\sqrt{c}}{p}.
$$
Computing the $D$-transform of $\mu_Z$  in (\ref{eq:mu z}) (see Example 3.1 in \cite{benaych2011svd} for the computation when $p=1$ from which the general $p$ answer can be easily deduced)  gives us the pertinent expression for $\woptbar{i}$ and $\overline{w}_{i}^{\rm eym}$ which match the expressions in Theorem \ref{conj:missing}. The $r = 1$ phase transition behavior for $\woptbar{1}$ follows from Theorem \ref{th:r1 soft thresholding}.

From the perturbation theory of singular values \cite[Theorem 3.3.16-(c), pp. 178]{hj91}, we have that
\begin{equation}\label{eq:lipschitz}
|\sigma_{i}(pS + Z +  \deltaS) - \sigma_{i}(pS + Z)| \leq \sigma_{1}(\deltaS),
\end{equation}
for $i = 1, \ldots, \min( m,n)$.  Consequently
$$|\weym{i} - \overline{w}_{i}^{\rm eym}| \leq \sigma_{1}(\deltaS),$$
so if we can show that $\sigma_{1}(\deltaS) \convas 0$ then we will have shown that $\weym{i} \convas \overline{w}_{i}^{\rm eym}$ and we have proved Theorem \ref{conj:missing}-a).

To prove that $\wopt{i} \convas \woptbar{i}$ we need a more involved argument that requires showing that we get the same limiting behavior when $Z + \deltaS$ is substituted for $Z$ in the bilinear forms on the left hand side of (\ref{eq:quad convergence}). We begin by noting that
\begin{align*}
 | u_{i}^{H} (w I_n - (Z+\deltaS)(Z+\deltaS)^{H})^{-1} u_{j}  &-u_{i}^{H} (w I_n - ZZ^{H})^{-1} u_{j}  |\\
  & \leq \sigma_{1}((w I_n - ZZ^{H})^{-1} -  (w I_n - (Z+\deltaS)(Z+\deltaS)^{H})^{-1})
\end{align*}
as a consequence of the variational characterization of the largest singular value. To make further progress, we shall utilize the resolvent identity\footnote{This identity can be verified by multiplying by $(wI-B)$ on the left and $(wI-A)$ on the right of the expressions on either side of the equality.} which states that
$$(wI - B)^{-1} - (wI - A)^{-1} = (wI - B)^{-1} (B-A)(wI-A)^{-1},$$
where $\Im w > 0$ and $A$ and $B$ are Hermitian matrices. Applying this identity with $A=ZZ^H$ and $B = (Z+\deltaS)(Z+\deltaS)^H$ yields
\begin{align*}
 | u_{i}^{H} (w I_n &- (Z+\deltaS)(Z+\deltaS)^{H})^{-1} u_{j}  -u_{i}^{H} (w I_n - ZZ^{H})^{-1} u_{j}  | \\
 & \leq \sigma_{1}((wI_n - (Z+\deltaS)(Z+\deltaS)^{H})^{-1} (\deltaS \deltaS^{H} + \deltaS Z^H + Z \deltaS^{H})(wI_n - ZZ^{H})^{-1} ) \\
&  \leq \dfrac{1}{|\Im w|^{2}} \cdot \sigma_{1}(\deltaS \deltaS^{H} + \deltaS Z^H + Z \deltaS^{H}).
\end{align*}
Since $\sigma_{1}(AB) \leq \sigma_{1}(A) \cdot \sigma_{1}(B)$  \cite[Theorem 3.3.16-(d), pp. 178]{hj91} and $\sigma_{1}(A+B) \leq \sigma_{1}(A) + \sigma_{1}(B)$  \cite[Theorem 3.3.16-(a), pp. 178]{hj91}, we have that
\begin{equation}
\sigma_{1}(\deltaS \deltaS^{H} + \deltaS Z^H + Z \deltaS^{H})  \leq \sigma_{1}^{2} (\deltaS) + 2 \sigma_{1}(Z) \sigma_{1}( \deltaS)  \leq 3 \sigma_{1}(Z) \sigma_{1}(\deltaS),
\end{equation}
if $\sigma_{1}(\deltaS) \leq \sigma_{1}(Z)$ thus leading to the inequality
\begin{equation}\label{eq:quad form}
 | u_{i}^{H} (w I_n - ZZ^{H})^{-1} u_{j} -  u_{i}^{H} (w I_n - (Z+\deltaS)(Z+\deltaS)^{H})^{-1} u_{j}|
  \leq \dfrac{3 \cdot \sigma_{1}(Z)}{|\Im w|^2} \sigma_{1}(\deltaS) .
 \end{equation}
Since $\sigma_1(Z) \convas b = \sqrt{p}\,(1+\sqrt{c}) < \infty$, if we can show that $\sigma_1(\deltaS) \convas 0$ we will have shown that the bilinear forms involving $u_i$ and $u_j$ exhibits the same limiting behavior as though $Z$ had i.i.d. Gaussian entries with zero mean and variance $p/m$. Repeating the argument would give us the analogous statement for the bilinear forms involving $v_i$ and $v_j$. To prove that $\sigma_{1}(\deltaS) \convas 0$, we first characterize $\mathbb{E}[\sigma_{1}(\deltaS)]$. From a theorem by Lata{\l}a \cite{latala2005some}, we have that
$$\mathbb{E}[\sigma_{1}(\deltaS)]  \leq C \left(  \max_{i} \sqrt{ \sum_{j} \mathbb{E}[\Delta S_{ij}^{2}]}  + \max_{j} \sqrt{ \sum_{i} \mathbb{E}[\Delta S_{ij}^{2}]}   +\sqrt[4]{ \sum_{ij} \mathbb{E}[\Delta S_{ij}^{4}]}\right) ,$$
where $C$ is a universal constant (that does not depend on $n$ or $m$). This gives us
\begin{equation}\label{eq:Esigma1DeltaS}
\mathbb{E}[\sigma_{1}(\deltaS)] \leq O\left(\dfrac{\textrm{$\log n$ factors}}{\sqrt{n}} \right).
\end{equation}
We note that
$$ |S_{ij}| \leq O\left(\dfrac{\textrm{$\log n$ factors}}{n} \right),$$
while
\begin{equation}\label{eq:max absSij}
\max_{i,j} |\Delta S_{ij}| \leq O\left(\dfrac{\textrm{$\log n$ factors}}{n}\right) =:K.
\end{equation}
Plugging in $i = 1$ in (\ref{eq:lipschitz}) we have
$$|\sigma_{1}(pS + Z +  \deltaS) - \sigma_{1}(pS + Z)| \leq 1 \cdot \sigma_{1}(\deltaS),$$
which implies that the largest singular value of a matrix is a $1$-Lipschitz function of the $n m $ entries of the matrix. Moreover, $\sigma_{1}(t\,A+(1-t)\,B) \leq t \sigma_{1}(A) + (1-t) \sigma_{1}(B) $, implying that the largest singular value is a convex, $1$-Lipschitz function. Since, by (\ref{eq:max absSij}), the entries of the $\deltaS$ are bounded, independent random variables, we can apply Talagrand's concentration inequality (see \cite[Theorem 2.1.13, pp. 73]{tao2012topics}) to obtain the tail bound
\begin{equation}\label{eq:tailSigma1DeltaS}
\textrm{Prob}\left( |\sigma_{1}(\deltaS) - \mathbb{E}[\sigma_{1}(\deltaS)]| > \epsilon \right) \leq 2 \exp\left(-c \, \frac{\epsilon^{2}}{K^{2}}\right) = 2 \exp\left(-c \,\frac{\epsilon^{2} n^{2}}{\textrm{$\log n$ factors}}\right).
\end{equation}
From  (\ref{eq:Esigma1DeltaS}), we have that $\mathbb{E}[\sigma_{1}(\deltaS)] \to 0$ as $n \to \infty$.  Moreover, the right-hand side of (\ref{eq:tailSigma1DeltaS}) is absolutely summable, \textit{i.e.},
$$ \sum_{n} 2 \exp\left(-c \,\frac{\epsilon^{2} n^{2}}{\textrm{$\log n$ factors}}\right) < \infty,$$
which implies, via the Borel-Cantelli lemma, that
\begin{equation} \label{eq:sigma1convas0}
\sigma_{1}(\deltaS) \convas 0
\end{equation}
Applying (\ref{eq:sigma1convas0}) to (\ref{eq:lipschitz}) yields the result that
$$\weym{i} = \sigma_{i}(pS + Z +  \deltaS)  \convas \sigma_{i}(pS + Z) = \overline{w}_{i}^{\rm eym}.$$
This proves Theorem \ref{conj:missing}-a). Moreover, from (\ref{eq:quad form}), we have that
$$ u_{i}^{H} (w I_n - (Z+\deltaS)(Z+\deltaS)^{H})^{-1} u_{j} \convas u_{i}^{H} (w I_n - ZZ^{H})^{-1} u_{j},$$
and by repeating the same argument we can show that
$$v_{i}^{H} (w I_m - (Z+\deltaS)^{H}(Z+\deltaS))^{-1} v_{j} \convas v_{i}^{H} (w I_m - Z^HZ)^{-1} v_{j}.$$
Using the same argument it can be shown that
$$ u_{i}^{H} (w I_n - (Z+\deltaS)(Z+\deltaS)^{H})^{-1} (Z+\deltaS) v_{j} \convas 0,$$
and
$$ v_{i}^{H} (w I_m - (Z+\deltaS)^H(Z+\deltaS))^{-1}  (Z+\deltaS)^H  u_{j} \convas 0.$$
Following the proofs in \cite{benaych2011svd}, the convergence of these bilinear forms implies that the almost sure limits of $\sigma_{i}(\wtX)$ and $\uhatu{i}{j}$ and $\vvhat{j}{i}$ for $i,j =1, \ldots, r$ are identical to the almost sure limits of $\sigma_{i}(\overline{X})$ and $\ubaru{i}{j}$ and $\vvbar{j}{i}$ for $i,j =1, \ldots, r$. Consequently, $\wopt{i} \convas \woptbar{i}$ and we have proved Theorem \ref{conj:missing}-b) and c).

\section{Justification for assumptions in Conjectures  \ref{conj:partial delocalization} and \ref{conj:full delocalization}}

A key aspect (see \cite[Lemma 4.1]{benaych2011svd}) in rigorously proving Conjectures \ref{conj:partial delocalization} and \ref{conj:full delocalization}
is understanding the behavior of expressions of the form
$$ u_{i}^{H} (z_j^2 I_n - XX^{H})^{-2} u_{i},$$
where $z_j$ is a singular value of $\wtX$ but not of $X$. Let $X  = U \Sigma V^H$ and $w  = U^H u_i$. Then
$$ u_{i}^{H} (z_j^2 I_n - XX^{H})^{-2} u_{i} = \sum_i \dfrac{|w_{i}|^2}{(z_j^2 - \sigma^2_i(XX^H))^2} \geq \dfrac{|w_j|^2}{(\sigma^2_{i+r}(XX^H)- \sigma^2_{i}(XX^H))^2}.$$
When $X$ has isotropically random singular vectors, $w_j = O(1/n)$ with high probability so if $z_j \in [a,b]$ and $\max_i \sigma_i(XX^H) - \sigma_{i+1}(XX^H)$ is bounded with probability by $O(\log n/n)$ in the bulk and the right hand side of the above expression will get unbounded (with $n$) resulting delocalization of the associated singular vectors. When $\mu_X$ exhibits a square root decay at the edge, then we expect the singular values at the edge to be spaced $O(n^{-2/3})$ apart with high probability so we might delocalization via the same argument. See \cite{nadakuditi2013most} for an exposition of some of these issues and
\cite{benarous13,BEKYY13} for recent results on the fine details of the spacing distribution of Wigner and Wishart random matrices.

%\bibliographystyle{plain}
%\bibliography{Asilomar_bib,spiked_bib}

\begin{thebibliography}{10}

\bibitem{bai2010spectral}
Z~Bai and J.~W.~Silverstein.
\newblock Spectral analysis of large dimensional random matrices,  2010.

\bibitem{baik2005phase}
J.~Baik, G.~Ben~Arous, and S.~P{\'e}ch{\'e}.
\newblock Phase transition of the largest eigenvalue for nonnull complex sample
  covariance matrices.
\newblock {\em The Annals of Probability}, 33(5):1643--1697, 2005.

\bibitem{baik2006eigenvalues}
J.~Baik and J.~W.~Silverstein.
\newblock Eigenvalues of large sample covariance matrices of spiked population
  models.
\newblock {\em Journal of Multivariate Analysis}, 97(6):1382--1408, 2006.

\bibitem{benarous13}
G.~Ben~Arous and P.~Bourgade.
\newblock Extreme gaps between eigenvalues of random matrices.
\newblock {\em The Annals of Probability}, vol. 41, no. 4, pp. 2648--2681, 2013.

\bibitem{benaych2011svd}
F. Benaych-Georges and R.~R.~Nadakuditi.
\newblock The singular values and vectors of low rank perturbations of large
  rectangular random matrices.
\newblock {\em Journal of Multivariate Analysis}, vol. 111, pp. 120--135, 2012.

\bibitem{b09}
F.~Benaych-Georges.
\newblock Rectangular random matrices, related convolution.
\newblock \emph{Probab. Theory Related Fields}, 144(3-4):471--515, 2009.

\bibitem{birnbaum2012minimax}
A. Birnbaum, I.~M.~Johnstone, B. Nadler, and D. Paul.
\newblock Minimax bounds for sparse pca with noisy high-dimensional data.
\newblock {\em arXiv:1203.0967}, 2012.

\bibitem{BEKYY13}
A.~Bloemendal, L.~Erd{\H{o}}s, A.~Knowles, H.-T. Yau, and J.~Yin.
\newblock Isotropic local laws for sample covariance and generalized wigner
  matrices.
\newblock arXiv: 1308.5729.

\bibitem{boucheron2013concentration}
 S.~Boucheron, G.~Lugosi and P.~Massart.
\newblock {C}oncentration {I}nequalities: {A} {N}onasymptotic {T}heory of {I}ndependence, Oxford University Press, 2013.


\bibitem{boutsidis2008svd}
C. Boutsidis and E. Gallopoulos.
\newblock {SVD} based initialization: A head start for nonnegative matrix
  factorization.
\newblock {\em Pattern Recognition}, 41(4):1350--1362, 2008.

\bibitem{breiman1995better}
L. Breiman.
\newblock Better subset regression using the nonnegative garrote.
\newblock {\em Technometrics}, 37(4):373--384, 1995.

\bibitem{cacciapuoti2012local}
C. Cacciapuoti, A. Maltsev, and B. Schlein.
\newblock Local {M}archenko-{P}astur law at the hard edge of sample covariance
  matrices.
\newblock {\em arXiv:1206.1730}, 2012.

\bibitem{cadzow1991enhanced}
J.~A.~Cadzow and D.~M.~Wilkes.
\newblock Enhanced rational signal modeling.
\newblock {\em Signal processing}, 25(2):171--188, 1991.

\bibitem{cai2010singular}
J.~F.~Cai, E.~J.~Cand{\`e}s, and Z.~Shen.
\newblock A singular value thresholding algorithm for matrix completion.
\newblock {\em SIAM Journal of Optimization}, 20(4):1956--1982, 2010.

\bibitem{Candes:2011:RPC:1970392.1970395}
E.~J.~Cand\`{e}s, X.~Li, Y.~Ma, and J.~Wright.
\newblock Robust principal component analysis?
\newblock {\em J. ACM}, 58(3):11:1--11:37, June 2011.

\bibitem{candes2010matrix}
E.~J.~Candes and Y.~Plan.
\newblock Matrix completion with noise.
\newblock {\em Proceedings of the IEEE}, 98(6):925--936, 2010.

\bibitem{candes2009exact}
E.~J.~Cand{\`e}s and B.~Recht.
\newblock Exact matrix completion via convex optimization.
\newblock {\em Foundations of Computational <athematics}, 9(6):717--772, 2009.

\bibitem{candes2010power}
E.~J.~Cand{\`e}s and T.~Tao.
\newblock The power of convex relaxation: Near-optimal matrix completion.
\newblock {\em IEEE Trans. on Information Theory}, 56(5):2053--2080,
  2010.

\bibitem{chandrasekaran2009sparse}
V.~Chandrasekaran, S.~Sanghavi, P.~A.~Parrilo, and A.~S.~Willsky.
\newblock Sparse and low-rank matrix decompositions.
\newblock In {\em Proc.  47th Annual Allerton Conference on Communication, Control, and Computing}, pages 962--967. IEEE, 2009.

\bibitem{chandrasekaran2011rank}
V.~Chandrasekaran, S.~Sanghavi, P.~A.Parrilo, and A.~S.~Willsky.
\newblock Rank-sparsity incoherence for matrix decomposition.
\newblock {\em SIAM Journal on Optimization}, 21(2):572--596, 2011.

\bibitem{chatterjee2012matrix}
S.~Chatterjee.
\newblock Matrix estimation by universal singular value thresholding.
\newblock {\em arXiv:1212.1247}, 2012.

\bibitem{chen2004recovering}
P.~Chen and D.~Suter.
\newblock Recovering the missing components in a large noisy low-rank matrix:
  Application to sfm.
\newblock {\em IEEE Trans. on Pattern Analysis and Machine Intelligence}, 26(8):1051--1063, 2004.

\bibitem{chu2004optimality}
M.~T.~Chu, F.~Diele, R.~Plemmons, and S.~Ragni.
\newblock Optimality, computation, and interpretation of nonnegative matrix
  factorizations.
\newblock In {\em SIAM Journal on Matrix Analysis}. Citeseer, 2004.

\bibitem{chu2003structured}
M.~T.~Chu, R.~E.~Funderlic, and R.~J.~Plemmons.
\newblock Structured low rank approximation.
\newblock {\em Linear algebra and its applications}, 366:157--172, 2003.

\bibitem{combettes92}
P.~L.~Combettes and J.~W.~Silverstein.
\newblock Signal detection via spectral theory of large dimensional random matrices.
\newblock {\em IEEE Trans. on Sig. Proc.}, vol. 8(40), pp. 2100--2105, 1992.

\bibitem{dasgupta1999learning}
S.~Dasgupta.
\newblock Learning mixtures of gaussians.
\newblock In {\em Foundations of Computer Science, 1999. 40th Annual Symposium
  on}, pages 634--644. IEEE, 1999.

\bibitem{d2008optimal}
A.~d'Aspremont, F.~Bach, and L.~El Ghaoui.
\newblock Optimal solutions for sparse principal component analysis.
\newblock {\em The Journal of Machine Learning Research}, 9:1269--1294, 2008.

\bibitem{deledalle2012risk}
C.-A.~Deledalle, S.~Vaiter, G.~Peyr{\'e}, J.~Fadili, and
  C.~Dossal.
\newblock Risk estimation for matrix recovery with spectral regularization.
\newblock {\em arXiv:1205.1482}, 2012.

\bibitem{drineas2006fast}
P.~Drineas, R.~Kannan, and M.~W.~Mahoney.
\newblock Fast monte carlo algorithms for matrices ii: Computing a low-rank
  approximation to a matrix.
\newblock {\em SIAM Journal on Computing}, 36(1):158--183, 2006.

\bibitem{eckart1936approximation}
C.~Eckart and G.~Young.
\newblock The approximation of one matrix by another of lower rank.
\newblock {\em Psychometrika}, 1(3):211--218, 1936.

\bibitem{el2007tracy}
N.~El~Karoui.
\newblock Tracy--widom limit for the largest eigenvalue of a large class of
  complex sample covariance matrices.
\newblock {\em The Annals of Probability}, 35(2):663--714, 2007.

\bibitem{erdHos2012universality}
L.~Erd{\H{o}}s and H.-T. Yau.
\newblock Universality of local spectral statistics of random matrices.
\newblock {\em Bull. Amer. Math. Soc}, 49:377--414, 2012.

\bibitem{fazel2008compressed}
M.~Fazel, E.~j.~Candes, B.~Recht, and P.~Parrilo.
\newblock Compressed sensing and robust recovery of low rank matrices.
\newblock In {\em 42nd {A}silomar {C}onference on {S}ignals, {S}ystems and {C}omputers, 2008 }, pp. 1043--1047. IEEE, 2008.

\bibitem{feral2009largest}
D.~F{\'e}ral and S.~P{\'e}ch{\'e}.
\newblock The largest eigenvalues of sample covariance matrices for a spiked
  population: diagonal case.
\newblock {\em Journal of Mathematical Physics}, 50:073302, 2009.

\bibitem{golub1987generalization}
G.~H.~Golub, A.~Hoffman, and G.~W.~Stewart.
\newblock A generalization of the {E}ckart-{Y}oung-{M}irsky matrix approximation
  theorem.
\newblock {\em Linear Algebra and Its Applications}, 88:317--327, 1987.

\bibitem{hachem2012subspace}
W.~Hachem, P.~Loubaton, X.~Mestre, J.~Najim, and P.~Vallet.
\newblock A subspace estimator for fixed rank perturbations of large random
  matrices.
\newblock {\em Journal of Multivariate Analysis}, 2012.

\bibitem{hj91}
R.~A.~Horn and C.~R.~Johnson.
\newblock {\em Topics in matrix analysis}.
\newblock Cambridge University Press, Cambridge, 1991.

\bibitem{hsu2012learning}
D.~Hsu and S.~M. Kakade.
\newblock Learning gaussian mixture models: {M}oment methods and spectral
  decompositions.
\newblock {\em arXiv:1206.5766}, 2012.

\bibitem{jenatton2010structured}
R.~Jenatton, G.~Obozinski, and F.~Bach.
\newblock Structured sparse principal component analysis.
\newblock In {\em International Conference on Artificial Intelligence and
  Statistics (AISTATS)}, 2010.

\bibitem{johnstone2009consistency}
I.~M.~Johnstone and A.~Lu.
\newblock On consistency and sparsity for principal components analysis in high
  dimensions.
\newblock {\em Journal of the American Statistical Association}, 104(486),
  2009.

\bibitem{johnstone2001distribution}
I.~M.~Johnstone.
\newblock On the distribution of the largest eigenvalue in principal components
  analysis.
\newblock {\em The Ann. of Statistics}, 29(2):295--327, 2001.

\bibitem{johnstone2006high}
I.~M.~Johnstone.
\newblock High dimensional statistical inference and random matrices.
\newblock In {\em Proceedings of the International Congress of Mathematicians:
  Madrid}, pages 307--333, 2006.

\bibitem{jolliffe2002principal}
I.~T.~Jolliffe.
\newblock {\em Principal component analysis}, volume~2.
\newblock Wiley Online Library, 2002.

\bibitem{kakade2012regularization}
S.~M.~Kakade, S.~Shalev-Shwartz, and A.~Tewari.
\newblock Regularization techniques for learning with matrices.
\newblock {\em The Journal of Machine Learning Research}, 98888:1865--1890,
  2012.

\bibitem{kannan2005spectral}
R.~Kannan, H.~Salmasian, and S.~Vempala.
\newblock The spectral method for general mixture models.
\newblock {\em Learning Theory}, pages 155--199, 2005.

\bibitem{kannan2009spectral}
R.~Kannan and S.~Vempala.
\newblock {\em Spectral algorithms}.
\newblock Now Publishers Inc, 2009.

\bibitem{karger2011budget}
D.~R.~Karger, S.~Oh, and D.~Shah.
\newblock Budget-optimal crowdsourcing using low-rank matrix approximations.
\newblock In {\em Communication, Control, and Computing (Allerton), 2011 49th
  Annual Allerton Conference on}, pages 284--291. IEEE, 2011.

\bibitem{keshavan2010matrix}
R.~H.~Keshavan, A.~Montanari, and S.~Oh.
\newblock Matrix completion from a few entries.
\newblock {\em , IEEE Trans. on Information Theory}, 56(6):2980--2998,
  2010.

\bibitem{klema1980singular}
V.~C.~Klema and A.~Laub.
\newblock The singular value decomposition: Its computation and some
  applications.
\newblock {\em IEEE Trans. on Automatic Control}, 25(2):164--176, 1980.

\bibitem{koltchinskii2011nuclear}
V.~ Koltchinskii, K.~Lounici, and A.~B.~Tsybakov.
\newblock Nuclear-norm penalization and optimal rates for noisy low-rank matrix
  completion.
\newblock {\em The Annals of Statistics}, 39(5):2302--2329, 2011.

\bibitem{kritchman2008determining}
S.~Kritchman and B.~Nadler.
\newblock Determining the number of components in a factor model from limited
  noisy data.
\newblock {\em Chemometrics and Intelligent Laboratory Systems}, 94(1):19--32,
  2008.

\bibitem{kritchman2009non}
S.~Kritchman and B.~Nadler.
\newblock Non-parametric detection of the number of signals: hypothesis testing
  and random matrix theory.
\newblock {\em IEEE Trans. on Signal Processing}, 57(10):3930--3941,
  2009.

\bibitem{langville2006initializations}
A.~N.~Langville, C.~D.~Meyer, R.~Albright, J.~Cox, and D.~Duling.
\newblock Initializations for the nonnegative matrix factorization.
\newblock In {\em Proceedings of the Twelfth ACM SIGKDD International
  Conference on Knowledge Discovery and Data Mining}, pages 23--26. Citeseer,
  2006.

\bibitem{latala2005some}
R.~Lata{\l}a.
\newblock Some estimates of norms of random matrices.
\newblock {\em Proc. of the American Math. Soc.}, vol. 133, no. 5, pp. 1273--1282, 2005.
	

\bibitem{le1960locally}
L.~M.~Le~Cam
\newblock  Locally asymptotically normal families of distributions: certain approximations to families of distributions and their use in the theory of estimation and testing hypotheses.
\newblock In {\em University of California Press}, vol. 3, no. 2, 1960.


\bibitem{li1997parameter}
Y.~Li, K.~J.~R.~Liu, and J.~Razavilar.
\newblock A parameter estimation scheme for damped sinusoidal signals based on
  low-rank {H}ankel approximation.
\newblock {\em IEEE Trans. on Signal Processing}, 45(2):481--486, 1997.

\bibitem{marchenko1967distribution}
V.~A.~Marchenko and L.~A.~Pastur.
\newblock Distribution of eigenvalues for some sets of random matrices.
\newblock {\em Matematicheskii Sbornik}, 114(4):507--536, 1967.

\bibitem{markovsky2008structured}
I.~Markovsky.
\newblock Structured low-rank approximation and its applications.
\newblock {\em Automatica}, 44(4):891--909, 2008.

\bibitem{mirsky1960symmetric}
L.~Mirsky.
\newblock Symmetric gauge functions and unitarily invariant norms.
\newblock {\em The quarterly journal of mathematics}, 11(1):50, 1960.

\bibitem{nadakuditi2011exploiting}
R.~R.~Nadakuditi.
\newblock Exploiting random matrix theory to improve noisy low-rank matrix
  approximation.
\newblock In {\em Signals, Systems and Computers (ASILOMAR), 2011 Conference
  Record of the Forty Fifth Asilomar Conference on}, pages 769--773. IEEE,
  2011.

\bibitem{nadakuditi2013most}
R.~R.~Nadakuditi.
\newblock When are the most informative components for inference also the
  principal components?
\newblock {\em arXiv:1302.1232}, 2013.

\bibitem{nadakuditi2012graph}
R.~R.~Nadakuditi and M.~E.~J.~Newman.
\newblock Graph spectra and the detectability of community structure in
  networks.
\newblock {\em Physical Review Letters}, 108(18):188701, 2012.

\bibitem{nadakuditi2013spectra}
R.~R.~Nadakuditi and M.~E.~J.~Newman.
\newblock Spectra of random graphs with arbitrary expected degrees.
\newblock {\em Physical Review E}, 87(1):012803, 2013.

\bibitem{nadakuditi2008sample}
R.~R.~Nadakuditi and A.~Edelman.
\newblock Sample eigenvalue based detection of high-dimensional signals in
  white noise using relatively few samples.
\newblock {\em IEEE Trans. on Signal Processing}, 56(7):2625--2638,
  2008.

\bibitem{nadler2010nonparametric}
B.~Nadler.
\newblock Nonparametric detection of signals by information theoretic criteria:
  performance analysis and an improved estimator.
\newblock {\em IEEE Trans. on Signal Processing}, 58(5):2746--2756,
  2010.

\bibitem{srebro2003weighted}
N.~Srebro and T.~Jaakkola.
\newblock Weighted low-rank approximations.
\newblock In {\em In 20th International Conference on Machine Learning},
  volume~20, page 720, 2003.

\bibitem{onatski09}
A.~Onatski.
\newblock Testing hypotheses about the numbers of factors in large factor models.
\newblock {\em Econometrica}, 77(5):1447--1479, 2009.

\bibitem{onatski2010determining}
A.~Onatski.
\newblock Determining the number of factors from empirical distribution of
  eigenvalues.
\newblock {\em The Review of Economics and Statistics}, 92(4):1004--1016, 2010.

\bibitem{passemier12}
D.~Passemier and J.-F.~Yao.
\newblock On determining the number of spikes in a high-dimensional spiked population model. 
\newblock {\em {R}andom {M}atrices: {T}heory and {A}pplications}, 1(1):1150002, 2012.

\bibitem{oymak2011finding}
S.~Oymak and B.~Hassibi.
\newblock Finding dense clusters via ``low rank+ sparse'' decomposition.
\newblock {\em arXiv:1104.5186}, 2011.

\bibitem{paul2007asymptotics}
D.~Paul.
\newblock {Asymptotics of sample eigenstructure for a large dimensional spiked
  covariance model}.
\newblock {\em Statistica Sinica}, 17(4):1617, 2007.

\bibitem{peche2009universality}
S.~P{\'e}ch{\'e}.
\newblock Universality results for the largest eigenvalues of some sample
  covariance matrix ensembles.
\newblock {\em Probability Theory and Related Fields}, 143(3-4):481--516, 2009.

\bibitem{pillai2011universality}
N.~S.~Pillai and J.~Yin.
\newblock Universality of covariance matrices.
\newblock {\em arXiv preprint arXiv:1110.2501}, 2011.

\bibitem{rohde2011estimation}
A.~Rohde and A.~B.~Tsybakov.
\newblock Estimation of high-dimensional low-rank matrices.
\newblock {\em The Annals of Statistics}, 39(2):887--930, 2011.

\bibitem{sanjeev2001learning}
A.~Sanjeev and R.~Kannan.
\newblock Learning mixtures of arbitrary gaussians.
\newblock In {\em Proceedings of the thirty-third annual ACM symposium on
  Theory of computing}, pages 247--257. ACM, 2001.

\bibitem{saunderson2012diagonal}
J.~Saunderson, V.~Chandrasekaran, P.~A.~Parrilo, and A.~S.~Willsky.
\newblock Diagonal and low-rank matrix decompositions, correlation matrices,
  and ellipsoid fitting.
\newblock {\em SIAM Journal on Matrix Analysis and Applications},
  33(4):1395--1416, 2012.

\bibitem{scharf1991svd}
L.~L.~Scharf.
\newblock The svd and reduced rank signal processing.
\newblock {\em Signal Processing}, 25(2):113--133, 1991.

\bibitem{shabalin2013reconstruction}
A.~.A.~Shabalin and A.~B.~Nobel.
\newblock Reconstruction of a low-rank matrix in the presence of {G}aussian noise.
\newblock{\em Journal of Multivariate Analysis}, 2013.

\bibitem{silverstein1995analysis}
J.~W.~Silverstein and S.-I. Choi.
\newblock Analysis of the limiting spectral distribution of large dimensional random matrices.
\newblock {\rm Journal of Multivariate Analysis}, vol. 54, no. 2, pp. 295--309, 1995.

\bibitem{soshnikov2002note}
A.~Soshnikov.
\newblock A note on universality of the distribution of the largest eigenvalues
  in certain sample covariance matrices.
\newblock {\em Journal of Statistical Physics}, 108(5-6):1033--1056, 2002.

\bibitem{tao2011recovering}
M.~Tao and X.~Yuan.
\newblock Recovering low-rank and sparse components of matrices from incomplete
  and noisy observations.
\newblock {\em SIAM Journal on Optimization}, 21(1):57--81, 2011.

\bibitem{tao2012topics}
T.~Tao.
\newblock Topics in random matrix theory.
\newblock vol. 132, AMS, 2012.

\bibitem{tipping1999mixtures}
M.~E.~Tipping and C.~M.~Bishop.
\newblock Mixtures of probabilistic principal component analyzers.
\newblock {\em Neural computation}, 11(2):443--482, 1999.

\bibitem{tufts1993estimation}
D.~W.~Tufts and A.~A.~Shah.
\newblock Estimation of a signal waveform from noisy data using low-rank
  approximation to a data matrix.
\newblock {\em Signal Processing, IEEE Transactions on}, 41(4):1716--1721,
  1993.

\bibitem{ulfarsson08}
M.~O.~Ulfarsson and V.~Solo.
\newblock Dimension estimation in noisy {PCA} with {SURE} and random matrix theory.
\newblock {\em IEEE Trans. on Sig. Proc.}, vol. 56(12), pp. 5804--5816, 2008.

\bibitem{vempala2002spectral}
S.~Vempala and G.~Wang.
\newblock A spectral algorithm for learning mixtures of distributions.
\newblock In {\em Foundations of Computer Science, 2002. Proceedings. The 43rd
  Annual IEEE Symposium on}, pages 113--122. IEEE, 2002.

\bibitem{wikes1988iterated}
D.~M.~Wikes and M.~H.~Hayes.
\newblock Iterated toeplitz approximation of covariance matrices.
\newblock In {\em Acoustics, Speech, and Signal Processing, 1988. ICASSP-88.,
  1988 International Conference on}, pages 1663--1666. IEEE, 1988.

\bibitem{zhang2012sparse}
Y.~Zhang, A.~d’Aspremont, and L. El~Ghaoui.
\newblock Sparse {PCA}: Convex relaxations, algorithms and applications.
\newblock In {\em Handbook on Semidefinite, Conic and Polynomial Optimization},
  pages 915--940. Springer, 2012.

\bibitem{zhang2002low}
Z.~Zhang, H.~Zha, and H.~Simon.
\newblock Low-rank approximations with sparse factors {I}: Basic algorithms and
  error analysis.
\newblock {\em SIAM Journal on Matrix Analysis and Applications},
  23(3):706--727, 2002.

\bibitem{zou2006sparse}
H.~Zou, T.~Hastie, and R.~Tibshirani.
\newblock Sparse principal component analysis.
\newblock {\em Journal of computational and graphical statistics},
  15(2):265--286, 2006.

\end{thebibliography}

\end{document}